\documentclass[11pt,a4paper,twoside]{article}
\usepackage[UKenglish]{babel}
\usepackage[T1]{fontenc}
\usepackage[utf8x]{inputenc}
\usepackage{latexsym,amsfonts,amsmath,amsthm,amssymb,mathrsfs}
\usepackage{fullpage}
\usepackage{xcolor,bm, bbm}
\usepackage{paralist}
\usepackage{todonotes}
\usepackage{dsfont}
\usepackage{float}
\usepackage{caption}
\usepackage{xfrac}
\usepackage{comment}
\usepackage{footnote}
\usepackage[multiple]{footmisc}

\usepackage{graphicx}
\usepackage{tikz,ifthen}
\usepackage{pgfplots}
\pgfplotsset{compat=1.18} 
\usepgfplotslibrary{fillbetween}

\usepackage{mathtools}

\DeclareMathOperator{\R}{\mathbb{R}}
\DeclareMathOperator{\N}{\mathbb{N}}

\DeclareMathOperator{\T}{\mathbb{T}}

\DeclareMathOperator{\rad}{rad}

\def \N {\mathbb{N}}
\def \ans {a_n^{(s)}}
\def \anl {a_n^{(r)}}
\def \ams {a_m^{(s)}}
\def \aml {a_m^{(r)}}

\def \ms {m^{(s)}}

\def \A {\mathcal{A}}

\def \ds {d^{(s)}}
\def \dl  {d^{(r)}}

\def \ms {m^{(s)}}

\def \ms {\mid^*}

\newtheorem{thm}{Theorem}

\newtheorem{lem}{Lemma}[section]

\newtheorem{cor}[thm]{Corollary}

\newtheorem*{rem*}{Remark}
\newtheorem*{thm*}{Theorem}
\newtheorem*{prob}{Open Problem}

\title{\bf Diophantine approximation with restricted denominators, metric gaps, and GCD sums of rough numbers}
\medskip

\author{Manuel Hauke-Treuer and Marc Munsch}

\newcommand{\Addresses}{{
  \bigskip
  \footnotesize

  Manuel Hauke-Treuer, \textsc{Institute of Analysis and Number Theory, TU Graz},\par 
  \textit{E-mail address:} \texttt{hauke@math.tugraz.at}\vspace{5mm}

  Marc Munsch, \textsc{Institut Camille Jordan, Université Jean Monnet},\par
   \textit{E-mail address:} \texttt{marc.munsch@univ-st-etienne.fr}
  }}
\date{}

\begin{document}
\maketitle

\begin{abstract}
    We consider the question of approximating random irrationals by rationals with denominators restricted to arbitrary sequences of distinct integers.  In particular, we show that for every such sequence $(a_n)_n$, for Lebesgue almost every $\alpha$, $\lVert a_n\alpha \rVert \leq \frac{1}{n}$ for infinitely many $n$, and prove that this is sharp. For increasing integer sequences, this can be improved to $\lVert a_n\alpha \rVert \leq \frac{1}{n (\log \log n)^2}$. We also prove that, for almost every \(\alpha\) and every sufficiently large \(N\), $\lVert a_n\alpha \rVert \leq \frac{(\log N)^{2 + \varepsilon}}{N}, n \leq N$ for all but finitely many $N$, which also improves the best-known bound on metric minimal spacings. One technical key innovation is to show that normalized GCD sums supported on sets of sufficiently rough numbers are bounded, a result that might be of independent interest. 
\end{abstract}

\section{Introduction and statement of results}
\subsection{Diophantine approximation along subsequences}
The cornerstone of Diophantine approximation is Dirichlet's theorem, in both its asymptotic and uniform version: writing $\lVert \cdot \rVert$ for the distance to the nearest integer, we have for all $\alpha \in \mathbb{R}$ that 

\begin{equation}\label{asymp_Dir}
    \lVert n\alpha\rVert < \frac{1}{n} \text{ for infinitely many $n \in \mathbb{N}$},
\end{equation}
as well as the stronger formulation 

\begin{equation}\label{uniform_Dir}
    \forall N \geq 1 \quad \exists n \leq N: \lVert n\alpha\rVert < \frac{1}{N}.
\end{equation}
Both variants are sharp up to constants; the so-called ``badly approximable numbers'' satisfy $\inf_{n \geq 1} n \lVert n\alpha\rVert > 0$, whereas the uniform variant in \eqref{uniform_Dir} cannot be replaced by $o(1/N)$ for \textit{any} irrational.\\
More can be said when one asks for improvements upon \eqref{asymp_Dir} or \eqref{uniform_Dir} for \textit{generic} $\alpha$, i.e., one asks for the best approximation rate for Lebesgue almost every $\alpha$. The starting point of this metric theory in the asymptotic setting \eqref{asymp_Dir} is Khintchine's theorem: given a function $\psi: \N \to [0,1/2)$, and defining\footnote{Throughout the remainder of this paper, i.m. means ``infinitely many.''}
\[W(\psi) := \left\{\alpha \in [0,1): \lVert n\alpha \rVert \leq \psi(n) \text{ for i.m.  $n \in \N$} \right\},\]
we obtain for every monotonically decreasing $\psi: \N \to [0,1/2)$,
\begin{equation}
\lambda(W(\psi))= \begin{cases}
0 &\text { if }  \sum\limits_{n \in \N} \psi(n) < \infty,\\
    1 &\text{ if } \sum\limits_{n \in \N}\psi(n) = \infty,\label{khintchine}
\end{cases}
\end{equation}
where $\lambda$ denotes the $1$-dimensional Lebesgue measure.

For uniform approximation, replacing the $\frac{1}{N}$ in \eqref{uniform_Dir} by $\frac{1-\varepsilon}{N}$ for an arbitrary $\varepsilon > 0$ is not possible for almost every $\alpha$ -- this follows e.g. immediately from the Davenport--Schmidt Theorem \cite{DS}. 
Thus in particular, writing

\[W_u(\psi) := \{\alpha \in [0,1): \exists N_0 \in \N: \forall N \geq N_0  \,\exists n \leq N: \lVert n \alpha\rVert < \psi(N)\},\]
we have for $\psi_c(n) := \frac{c}{n}$ that\footnote{In fact, we note that a more recent work by Kleinbock--Wadleigh \cite[Theorem $1.8$]{KW}
provides an even more precise $0-1$-law for $W_u(\psi)$ under the assumptions of $\psi$ being decreasing, 
and $n \mapsto n\psi(n)$ increasing.}
\begin{equation}
\lambda(W_u(\psi_c))= \begin{cases}
0 &\text { if }  c < 1,\\
    1 &\text{ if } c \geq 1.\label{kw_khintchine}
\end{cases}
\end{equation}

The main purpose of the present article is to study how close we can get to  Khintchine-type laws such as \eqref{khintchine} and \eqref{kw_khintchine} when only allowing approximations with prescribed denominators $(a_n)_{n \in \mathbb{N}}$.
More formally, let $\psi: \N \to [0,1/2)$ be a decreasing function, and let $\A = (a_n)_{n \in \mathbb{N}}$ be a sequence of distinct integers. We define

\begin{equation}
    \label{W_normal}
W(\psi,\A) := \left\{\alpha \in [0,1): \lVert a_n\alpha \rVert \leq \psi(n) \text{ for i.m.\quad } n \in \mathbb{N} \right\},\end{equation}
as well as
\begin{equation}\label{W_unif}W_{u}(\psi,\A) := \left\{\alpha \in [0,1): \exists N_0\in \N: \forall N\geq N_0 \,\exists n \leq N:\lVert a_n\alpha \rVert \leq \psi(N)\right\}.\end{equation}
A na\"{\i}ve expectation might be that \eqref{khintchine} and \eqref{kw_khintchine} generalize to $W(\psi,\A)$ and $W_u(\psi,\A)$ in place of $W(\psi) = W(\psi,\mathbb{N})$ or $W_u(\psi) = W_u(\psi,\N)$ for arbitrary $\A$. By the first Borel--Cantelli lemma we obtain
\[\sum_{n \in \mathbb{N}}\psi(n) < \infty \implies \lambda(W(\psi,\A)) = 0.\]
Similarly, by a simple union bound and using Lebesgue density, one can show that for $c < \frac{1}{2}$, $\lambda(W_u(\psi_c,\A)) = 0.$ 
Thus we focus for the remainder of the article on criteria on $\psi$ that imply
$\lambda(W(\psi,\A)) = 1$ and $\lambda(W_u(\psi,\A)) = 1$, respectively, which is the far more delicate question. More precisely, we aim for the best possible rates for $\psi$ such that
\begin{equation}
    \label{adm_psi}
    \lambda(W(\psi,\A)) = 1 \text{ for every integer sequence $\A$,}
\end{equation}
respectively
\begin{equation}
    \label{adm_psi_unif}
    \lambda(W_u(\psi,\A)) = 1 \text{ for every integer sequence $\A$}.
\end{equation}

\subsection{Asymptotic approximation}
Let us focus first on asymptotic approximation \eqref{adm_psi}.
Note that by setting 
\begin{equation}\theta(q) := \begin{cases}
    \psi(n), &\text{ if } q = a_n,\\
    0 &\text{ if } q \notin \{a_n: n \in \N\},
\end{cases}
\label{eq:khintchine_transl}
\end{equation}
a Khintchine-type behaviour for $W(\psi,\A)$ is equivalent to asking whether Khintchine's theorem is true without monotonicity assumption (cf. Harman's monograph \cite[Section 2.6]{H98}). This was ruled out by the counterexample of Duffin--Schaeffer \cite{DS41}, 
hence \eqref{khintchine} cannot be generalized to \eqref{W_normal}.
In this article, we show that for arbitrary integer sequences $\A$, the best admissible approximation rate is $\psi(n) = \frac{1}{n}$. 

\begin{thm}\label{nonincr_thm}
    For every sequence of distinct integers $\A = (a_n)_n$, we have
    for almost every $\alpha$,
    \begin{equation}\label{eq:1nrate}\liminf_{n \to \infty} n \lVert a_n \alpha \rVert = 0.\end{equation}
    This is sharp in the following sense:
    for every function $f: \N \to [0,\infty)$ with $\lim_{n \to \infty} f(n) = \infty$, there exists 
    a sequence of distinct integers $\A = (a_n)_n$ such that for almost every $\alpha$, 
        \begin{equation}\label{eq:1n_sharp}\liminf_{n \to \infty} f(n) n \lVert a_n \alpha \rVert > 0.\end{equation}
\end{thm}

\begin{rem*}
    Theorem \ref{nonincr_thm} has, to the best knowledge of the authors, not been stated prior to this article. The first part \eqref{eq:1nrate} however follows from the proof method of \cite[Theorem 1]{AEM}, observing that assumptions on $\A$ made in their theorem statement are in fact never used in the proof. To provide a clean self-contained proof for \eqref{eq:1nrate}, we give, while relying partially on ideas from \cite{AEM}, a slightly different argument that is both shorter and closer to the strategy used to prove Theorem \ref{main_thm_asymp}, where both proof and statement appear to be novel. The best previously \textit{stated} bound was $\frac{n}{\exp\left(C \sqrt{\log n \log \log \log n/\log \log n}\right)}$ in place of the $n$ from \eqref{eq:1nrate}, see~\cite{BS_15}. Since it is in some sense hard to correctly attribute the currently best bounds to a precise article, we refer the reader to the discussion in Subsection \ref{sec:intro_gcd} below.
\end{rem*}

We remark that the sequences $\A$ that are constructed for \eqref{eq:1n_sharp} are far from being \textit{monotonically increasing}. The potentially more natural question is to consider \eqref{W_normal} for strictly \textit{increasing} integer sequences $\A$ only. Using the same translation \eqref{eq:khintchine_transl}, a Khintchine-type result for increasing sequences  $\A$ is equivalent to Khintchine's theorem being true for approximation functions $\theta$  that are monotonically decreasing on their support. The recent counterexamples in \cite[Theorem $1.2$]{CHPR25} show that even in this setting, an analogue of \eqref{khintchine} for $W(\psi,\A)$ fails.
Our main result concerning asymptotic approximation 
however shows that for strictly increasing sequences $\A$, we can indeed improve upon \eqref{eq:1nrate}.
\begin{thm}\label{main_thm_asymp}
Let $\A=(a_n)_{n\in\mathbb N}$ be a strictly increasing
sequence of positive integers.
Then, for Lebesgue almost every $\alpha\in\R$,
\[ \liminf_{n \to \infty} 
 n (\log \log n)^2\lVert a_n\alpha\rVert = 0.
\]
\end{thm}

We will see in Section \ref{sec:intro_gcd} how the rate in Theorem \ref{main_thm_asymp} might be further improved by stronger bounds for certain GCD sum estimates.
We further note that a direct application of the Koukoulopoulos--Maynard theorem \cite{KM20} (formerly known as Duffin--Schaeffer conjecture) as well as  the classical bound\footnote{Throughout this article, $\varphi(n)$ stands for the Euler totient function.} $\frac{\varphi(a_n)}{a_n} \gg (\log \log a_n)^{-1}$ allows 
for better estimates, provided $(a_n)_n$ does not grow too quickly. Combined with Theorem \ref{main_thm_asymp}, this leads to the following statement.

\begin{cor}\label{asymp_thm_cases}
Let $\A=(a_n)_{n\in\mathbb N}$ be a strictly increasing
sequence of positive integers.
Then, for Lebesgue almost every $\alpha\in\R$, the inequality
\[
\lVert a_n\alpha\rVert
\leq \frac{1}{n}
\min\left\{
    \frac{1}{(\log\log n)^{2}},
    \frac{\log\log a_n}{\log n\,\log\log n}
\right\}
\] holds for infinitely many $n\in\mathbb N$.
    \end{cor}

\subsection{Uniform approximation and the metric theory of minimal spacings}\label{subsec:spacing}
We will now focus on the question of \textit{uniform} approximation, i.e., sets $W_u(\psi)$, and ask how far away we are from the rate $\psi(N) = \frac{1}{N}$, as stated in \eqref{kw_khintchine} for $\A = \N$. Our main theorem on uniform approximation shows that, even for not necessarily increasing integer sequences, we reach \eqref{kw_khintchine} up to logarithmic factors.

 \begin{thm}\label{main_thm_unif}
Let $\A=(a_n)_{n\in\mathbb N}$ be an arbitrary sequence of distinct positive integers, and let $\varepsilon>0$ be arbitrary.
Then, for Lebesgue almost every $\alpha\in\R$, for all
sufficiently large $N\in\mathbb N$, there exists $1\leq n\leq N$
such that
\[
 \lVert a_n\alpha\rVert
 \leq \frac{(\log N)^{2+\varepsilon}}{N}.
\] 
\end{thm} 
Theorem \ref{main_thm_unif} improves upon the previously best-established bound of 
\[\psi(N) = \frac{\exp\left( (2\sqrt{2} + o(1)) \sqrt{\log N \log \log \log N/\log \log N}\right)}{N},\]
where we again refer the reader to Section \ref{sec:intro_gcd}
for the historical developments. We note that, in contrast to Theorem \ref{nonincr_thm}, no paralleling result from \cite{AEM} exists.

Similarly to the asymptotic question, the Duffin--Schaeffer theory, in this case in its quantitative version due to Koukoulopoulos--Maynard--Yang \cite{KMY24}, allows again for improvements for slowly growing sequences.

\begin{cor}\label{unif_thm_cases}
Let $\A=(a_n)_{n\in\mathbb N}$ be a strictly increasing 
 sequence of integers, and let $\varepsilon>0$ be arbitrary.
Then, for Lebesgue almost every $\alpha\in\R$ and all
sufficiently large $N\in\mathbb N$, there exists $1\leq n\leq N$
such that
\[
 \lVert a_n\alpha\rVert
 \leq \frac{(\log N)^{1+\varepsilon}}{N}\min\left\{\log N, \log \log a_N\right\}.
\]
    \end{cor}
 
\begin{rem*}
Note that the effectiveness of the Duffin--Schaeffer approach in Corollaries \ref{asymp_thm_cases} and \ref{unif_thm_cases} reaches its limits at different scales, depending on how slowly the sequence grows.
If $(a_n)_n$ grows e.g. at most polynomially, then Corollary \ref{asymp_thm_cases} allows asymptotic approximation rates of quality $\frac{1}{n \log n}$, which misses the best possible bound (up to triple-logs) by at most a factor of $ \log \log n$. On the other hand, for the uniform question, for any sequence with growth rate $a_n = O(\exp\{\exp((\log n)^{o(1)})\})$, Corollary \ref{unif_thm_cases} allows to show uniform approximations by
$\frac{(\log N)^{1 + \varepsilon}}{N}$, and this rate does not improve under stronger growth assumptions.
\end{rem*}

Theorem \ref{main_thm_unif} has direct applications to the study of local statistics of sequences, in particular to the metric theory of minimal spacings (we refer to \cite{AEM} for more information and references). 
For $(a_n)_{n \in \mathbb{N}}$ a sequence of distinct integers, we write 
\[\mathcal{A}_N = \{a_n: 1 \leq n \leq N\} \qquad  \text{and} \qquad C_N = \#(\mathcal{A}_N - \mathcal{A}_N). \] Further for $\alpha \in \mathbb{R}$, we define
\[\delta_{\min}^{\alpha}(N) := \min\{\lVert (a_n-a_m)\alpha\rVert: 1 \leq n \neq m \leq N\} = 
\min\{\lVert c_n\alpha\rVert: c_n \in (\mathcal{A}_N - \mathcal{A}_N)\setminus\{0\}\}.\]
Theorem \ref{main_thm_unif} gives the following corollary.
\begin{cor}\label{cor_alln}
   Let $\A=(a_n)_{n\in\mathbb N}$ be a strictly increasing  sequence of integers, and let $\varepsilon >0$ be arbitrary. Then, for Lebesgue almost every $\alpha \in \mathbb{R}$ and all sufficiently large $N \in \mathbb{N}$,  we have
    \begin{equation}\label{uppbnd_uni} \delta_{\min}^{\alpha}(N) \leq \frac{(\log N)^{1 + \varepsilon}}{C_N}\min\left\{\log N, \log \log a_N\right\}. \end{equation}
\end{cor} This improves upon results of Rudnick \cite{Rudnick} 
 who obtained\footnote{In \cite{Rudnick}, the estimate takes a weaker form, but it was shown in \cite{AEM} how to adapt the approach to the stated threshold.} $\frac{\exp\left((2\sqrt{2}+o(1))\sqrt{\frac{\log N \log \log \log N}{\log \log N}}\right)}{C_N}$. 
 We note that the assumption of $(a_n)_{n}$ being increasing is only necessary for the case of \\$\min\{\log N, \log \log a_N\} = \log \log a_N$, since the inequality 
 $ \delta_{\min}^{\alpha}(N) \leq \frac{(\log N)^{2 + \varepsilon}}{C_N}$ follows immediately from Theorem \ref{main_thm_unif} for \textit{arbitrary} sequences of distinct integers. \\
 Corollary \ref{cor_alln} makes progress towards the bound conjectured in \cite[Equation (7)]{AEM} which postulates that we could replace the right hand side of \eqref{uppbnd_uni} by $\frac{(\log N)^{\varepsilon}}{C_N}$. 
 It remains possible that the approximation function in both Theorem \ref{main_thm_unif} as well as Corollary \ref{cor_alln} can be replaced by $\frac{(\log N)^{\varepsilon}}{N}$ and $\frac{(\log N)^{\varepsilon}}{C_N}$, respectively. However, our method of an $L^2$-approach reaches its limit at $(\log N)^{1 + \varepsilon}$. We believe that in order to go significantly lower, one would need to establish a new approach that potentially includes higher moments.
 
 We remark at this point that the proof of Corollary \ref{cor_alln} does not make use of the fact that the sequence $(c_n)_n$ arises as a difference sequence, but treats it as an arbitrary sequence of distinct integers. In the related ``infinitely many'' setup (cf. Theorem \ref{nonincr_thm}), the above-mentioned result of \cite[Theorem 1]{AEM} proves that 
 $\delta_{\min}^{\alpha}(N) \leq \frac{1}{C_N}$ infinitely often,
 again not exploiting the particular structure of $(c_n)_n$.
 As witnessed by the second part of Theorem \ref{nonincr_thm}, this cannot be improved without exploiting the structure of a difference sequence. 
 It remains open whether the additive structure of a difference set can be used to improve the bounds for minimal gaps. This question is related to the classical sum-product phenomenon in additive combinatorics (see e.g. \cite[Section 8.3]{TV06}). Broadly speaking, the heuristic is that a finite set of integers cannot be simultaneously additively and multiplicatively structured. Clearly, a difference set has an additive structure, whereas the counterexamples presented in Theorem \ref{nonincr_thm} as well as maximizing ensembles for GCD sums contain a lot of multiplicative structure. For concrete results in this direction, we refer to Shkredov \cite{S16} who proved that difference sets have large product sets, and Bloom and Walker \cite{bw20}, who used a related sum–product estimate to bound GCD sums for pair correlations. In their setting, each difference is weighted by its multiplicity, whereas for minimal gaps each distinct difference counts only once.
Thus it remains possible that the rates for $\delta_{\min}^{\alpha}(N)$ from Corollary \ref{cor_alln} and \cite[Theorem 1]{AEM} can be improved by such input from additive combinatorics.

 \subsection{Historical bounds on the problem and the theory of GCD sums}\label{sec:intro_gcd}

The problem of Diophantine approximation along arbitrary integer sequences is historically closely connected
to the theory of GCD sums (sometimes referred to as ``Gál sums'' due to Gál's seminal work \cite{Gal}). LeVeque's work \cite{L76} on the study of metric
uniform distribution and dilated functions provided the first connection to such sums. In particular, the results there established $\psi(n) = \frac{1}{n^{1/2-\varepsilon}}$ for \eqref{adm_psi}, which was later improved by Harman \cite{H86}. The work of Dyer--Harman \cite{DH86} and Harman \cite{H86}
(see also Chapter 3 in Harman's monograph \cite{H98}) shows how bounds for GCD sums
yield approximation functions $\psi$ that satisfy \eqref{adm_psi} and \eqref{adm_psi_unif}. More precisely, if
$G: \N \to \mathbb{R}$ is non-decreasing, $\lim_{N \to \infty} \frac{G(2N)}{G(N)} = 1$, and $G$ satisfies the bound
\begin{equation}
    \label{def_gn}
 \limsup_{N \to \infty} \frac{1}{N G(N)}\sup_{\substack{\mathcal{M} \subseteq \N\\ |\mathcal{M}| = N}}\sum_{m,n\in \mathcal{M}}
 \frac{\gcd(m,n)}{\sqrt{mn}} < \infty,
\end{equation} 
then for any decreasing $\psi: \mathbb{N} \to [0,1/2)$ such that

\begin{equation}\label{extra_div_G}
    \sum_{n = 1}^{+\infty} \frac{\psi(n)}{G(n)} = \infty,
\end{equation}
we have $\lambda(W(\psi,\A)) = 1$.

Dyer--Harman \cite{DH86} established that
$G(N)=\exp(C\log N/\log\log N)$ satisfies \eqref{def_gn}, which was improved by Aistleitner, Berkes and Seip \cite{ABS_15},
and subsequently by Bondarenko and Seip \cite{BS_15,BS17}.
Finally, de la Bret\`eche and Tenenbaum \cite{BT} determined
the precise leading constant in the exponent, proving that
\begin{equation}\label{gcd_asymp}
 \sup_{\substack{\mathcal{M} \subseteq \N\\ |\mathcal{M}| = N}}\frac{1}{N}\sum_{m,n\in \mathcal{M}}
 \frac{\gcd(m,n)}{\sqrt{mn}}
 =
 \exp\left(
 (2\sqrt{2}+o(1))
 \sqrt{\frac{\log N\log\log\log N}{\log\log N}}
 \right).
\end{equation}

Similarly, up to a logarithmic loss that does not affect the leading
exponential scale of these bounds, such argument can be applied for the question of uniform approximation. Using \eqref{gcd_asymp}, one may therefore derive that

\[\psi(N) = \exp\left(
 (2\sqrt{2}+o(1))
 \sqrt{\frac{\log N\log\log\log N}{\log\log N}}
 \right)/N\]
is admissible for \eqref{adm_psi_unif}, which was (to the best knowledge of the authors) the best-known bound prior to this article.

We pause here to remark that large values of GCD sums have been applied in various other areas of mathematics beyond the application in metric Diophantine approximation described here. In particular, there is a line of research that identifies maximal values of GCD sums as spectral norms of Gram matrices associated with systems of dilated functions in the Hilbert space $L^2(0,1)$; later, an equivalent interpretation in terms of multiplication operators on Hardy spaces of Dirichlet series was intensively studied. The foundational paper for this branch of research is due to Wintner \cite{Wintner44}, with important later contributions by Hedenmalm--Lindqvist--Seip \cite{HLS97} and Lindqvist--Seip \cite{LS98}. The connection with Dirichlet series already foreshadowed the appearance of GCD sums in extreme value problems for the Riemann zeta function and $L$-functions, in a form that was developed by Soundararajan \cite{sound} and Hilberdink \cite{Hilb}, first made fully explicit by Aistleitner \cite{Ais}, and later applied by, e.g., Bondarenko--Seip \cite{BS17} and de la Bretèche-Tenenbaum \cite{BT}.

Returning to the problem of metric Diophantine approximation, the GCD sum advances  described above occurred before Koukoulopoulos--Maynard \cite{KM20} proved the  Duffin--Schaeffer conjecture. As mentioned, while $\frac{\varphi(a_n)}{a_n} \gg (\log \log a_n)^{-1}$ allows for good bounds for reasonably slowly-growing sequences (see Corollaries \ref{asymp_thm_cases} and \ref{unif_thm_cases}), the term $\frac{\varphi(a_n)}{a_n}$ might become arbitrarily small in terms of $n$; therefore, the Koukoulopoulos--Maynard theorem itself gives no result whatsoever for \textit{arbitrary} $\A$. Nevertheless, the \textit{tools} developed in \cite{KM20} turned out useful for \eqref{adm_psi},  as used in 
\cite[Theorem 1]{AEM}: instead of using sets of (completely) reduced fractions as in the classical Duffin-Schaeffer setting, (larger) sets were considered where only fractions whose numerator and denominator share a common \textit{small} prime factor are removed (see Section \ref{sec:heuristic} for a detailed discussion on the proof strategy). In that way, the bound $\psi(n) = \frac{1}{n}$ can be obtained for \eqref{adm_psi}, following the proof of \cite[Theorem 1]{AEM}. For the more natural consideration of increasing integer sequences (cf. Theorem \ref{main_thm_asymp}), a closer examination and a refined approach lead to considering the following optimisation problem concerning GCD sums supported on rough numbers.\footnote{We say that a set $\mathcal{M} \subseteq \N$ is $T$-rough if $P^{-}(n)\geq T$ for any $n\in \mathcal{M}$ where $P^{-}(n)$ stands for the smallest prime factor of $n$.} 
Precisely, we aim for the minimal growth rate for a function 
 $T: \mathbb{N} \to \mathbb{R}$ (compare with \eqref{def_gn})
so that
\begin{equation}\label{GCD_rough_bound}
\limsup_{N \to \infty} \sup_{\substack{\mathcal{M} \subseteq \N,\,\,\,T(N)-\text{rough}\\|\mathcal{M}| = N}} \frac{1}{N}\sum_{n,m \in \mathcal{M}}\frac{\gcd(n,m)}{\sqrt{nm}} < \infty.
\end{equation}
With this defining inequality at hand, we can state a more flexible version of Theorem \ref{main_thm_asymp}, which shows that an improvement on the rates for $T$ that satisfy \eqref{GCD_rough_bound} (and a technical doubling assumption that is satisfied by all natural choices) immediately provides an improvement on the rate for the asymptotic approximation quality along increasing integer sequences.

\begin{thm}\label{gen_asymp}
Let $\A=(a_n)_{n\in\mathbb N}$ be a strictly increasing
sequence of positive integers, and let 
$T: \mathbb{N} \to (1,+\infty)$ be an increasing function satisfying \eqref{GCD_rough_bound}, as well as 
\begin{equation}\label{doubling}\limsup_{n \to \infty} \frac{\log T(n^2)}{\log T(n)} < \infty.\end{equation}
Then for any decreasing $\psi:\N\to[0,1/2)$ such that
\begin{equation}\label{extra_div_thm}
 \sum_{n=1}^{+\infty} \frac{\psi(n)}{ \log T(n)}=\infty,
\end{equation}  for Lebesgue almost every $\alpha\in\R$, we have
\[
 \lVert a_n\alpha\rVert
 \leq \psi(n)
\]
for infinitely many $n\in\mathbb N$.
\end{thm}

The condition \eqref{extra_div_thm} should be compared to \eqref{extra_div_G}. Note that already the trivially admissible choice $T(N) = N^2$ provides the (non-trivial) rate $\psi(n) = \frac{1}{n\log \log n}$, which already improves upon the best previously established rate. However, a closer examination of \eqref{GCD_rough_bound} allows us to improve upon the trivial bound on an exponential scale, a result that might be of independent interest in multiplicative number theory.

\begin{thm}\label{doublelog_thm}
   Let $T(N) = \exp\left(\frac{20\log N}{\log \log N}\right)$ where $N$ is a sufficiently large integer. Then for any $T(N)$-rough set $\mathcal{M} \subseteq \mathbb{N}$ with $|\mathcal{M}|=N$ and any non-negative weights $(c_n)_{n \in \mathcal{M}}$, we have
    \begin{equation}\label{eq:weight_bd}
    \sum_{n,m\in \mathcal{M}} \frac{\gcd(n,m)}{\sqrt{nm}}\min\{c_n,c_m\} \ll \sum_{m \in \mathcal{M}}c_m.
    \end{equation}
    In particular, we have
    \begin{equation}\label{eq:unw_bd}\frac{1}{N}\sum_{n,m \in \mathcal{M}}\frac{\gcd(n,m)}{\sqrt{nm}} \ll 1.\end{equation} 
\end{thm}
A combination of Theorem \ref{gen_asymp} and Theorem \ref{doublelog_thm}, that form the technical core of this article immediately proves Theorem \ref{main_thm_asymp}.
Indeed, choosing $T(n) = \exp(20 \log n/\log \log n)$ ensures by Theorem \ref{doublelog_thm} that \eqref{GCD_rough_bound} holds. Thus 
$\psi_{\varepsilon}(n) = \frac{\varepsilon}{n (\log \log n)^2}$ allows for an application of Theorem \ref{gen_asymp}, proving for increasing $\A$ that 
$\lambda(W(\psi_{\varepsilon},\A)) = 1$. Since the sets $(W(\psi_{\varepsilon},\A))_{\varepsilon > 0}$ are nested, the result follows by continuity of measures.
In the converse direction, we show that our method of proof could never replace the rate from Theorem \ref{main_thm_asymp} with $\psi(n) = \frac{1}{n ( \log n)^ {\varepsilon}}$. It remains unclear whether such a rate is admissible for arbitrary increasing sequences.

\begin{thm}\label{unbounded_GCD}
Let $\varepsilon > 0$ be arbitrary, $N \in \N$ and let $T(N) = \exp((\log(N))^{1-\varepsilon})$. Then there exists $\mathcal{M} \subseteq \mathbb{N}$ with $|\mathcal{M}| = N$ consisting only of $T(N)$-rough numbers such that for every $C > 0$, 
\[\sum_{n,m \in \mathcal{M}} \frac{\gcd(m,n)}{\sqrt{mn}} \gg_{C,\varepsilon} N (\log N)^C.\]
\end{thm}

While Theorems \ref{doublelog_thm} and \ref{unbounded_GCD} provide upper and lower bounds for $T$ satisfying \eqref{GCD_rough_bound}, a sharp rate remains an open problem, which again might be of independent interest.

\begin{prob}
    Find sharp growth conditions for $T: \N \to (1,+\infty)$ such that \eqref{GCD_rough_bound} holds.
\end{prob}

\begin{rem*}
We remark that we did not aim to optimize the bound for Theorem \ref{unbounded_GCD} too much. It might well be that the ensembles $\mathcal{M}$ provided in the proof allow for quantifying $\varepsilon \to 0$, giving a sharper lower bound for the above Open Problem.
Conversely, we remark that the method to find the upper bound provided in Theorem \ref{doublelog_thm} does not allow for an improvement (except for the constant $20$ which is not optimized either).
\end{rem*}

\begin{rem*}
We remark that the weighted and unweighted formulations \eqref{eq:weight_bd} and \eqref{eq:unw_bd} are in fact equivalent, see Lemma \ref{lem:weighted_GCD}. However, replacing \eqref{eq:weight_bd} with the symmetric
\begin{equation}\label{spectral_}\sum_{n,m \in \mathcal{M}}\frac{\gcd(n,m)}{\sqrt{nm}}c_n^{1/2} c_m^{1/2} \ll \sum_{n \in \mathcal{M}} c_n,\end{equation}
which is the more classically studied object due to its connections to spectral norms of Gram matrices described above, makes a difference here: for every $T = N^{o(1)}$, one can find a $T$-rough set $\mathcal{M}$ and weights $(c_i)_i \geq 0$ such that \eqref{spectral_} fails. This is a technical explanation for the  discrepancy of the rates admissible for increasing $\A$ (see Theorem \ref{main_thm_asymp}) in contrast to the sharp rate for arbitrary $\A$ (see Theorem \ref{nonincr_thm}). For details on that matter, we refer to the discussion at the end of Section \ref{subsec:asymp}.
\end{rem*}

\subsection*{Organization of the paper}

The remainder of the paper is organized as follows. In Section
\ref{sec:heuristic}, we present the proof strategy in context of previous methods, highlighting how our proof differs from previous approaches. 
In Section \ref{sec:GCD}, we prove the results on GCD sums (Theorems \ref{doublelog_thm} and \ref{unbounded_GCD}). In Section \ref{sec:asymptotic}, we prove the results on asymptotic approximation (Theorems \ref{nonincr_thm}, \ref{main_thm_asymp}, \ref{gen_asymp} and Corollary \ref{asymp_thm_cases}), whereas in Section \ref{sec:uniform}, we prove the results on uniform approximation (Theorem \ref{main_thm_unif} and Corollaries \ref{unif_thm_cases} and \ref{cor_alln}).

\section*{Acknowledgements}
M.H.-T. was funded by the Austrian Science Fund (FWF) project
10.55776/ESP5134624. Parts of this project were carried out while M.H.-T.
visited M.M. in November 2024. He thanks University Jean Monnet for their hospitality. M.H.-T.
further thanks Andrii Bondarenko for various enlightening discussions concerning GCD sums, which turned out to be particularly useful for the example provided in Theorem \ref{unbounded_GCD}. Further, he thanks Christoph Aistleitner for various useful comments on an earlier version of this manuscript. During the preparation of this work, M.M. was supported by AAP Recherche 2025 UJM ``Comportements al{\'e}atoires en arithm{\'e}tique''. Moreover, this research was funded in whole or in part by the French National Research Agency (ANR) under project number ANR-25-CE40-1961-01.

\subsection*{AI declaration}
We did not use AI for any of the mathematical work or the writing of this manuscript. It is solely the work of the human authors.

\section{Heuristics and the general strategy}\label{sec:heuristic}

In this section, we provide heuristics and some historical background for the three different aspects considered in this paper:

\begin{itemize}
    \item Ideas for asymptotic approximation,
    \item Ideas for uniform approximation,
    \item Ideas for rough GCD sums.
\end{itemize}

\subsection{Ideas for asymptotic approximation}\label{subsec:asymp}

Here, we review the historical background and the proof strategy for
Theorem \ref{gen_asymp}, and hence Theorem \ref{main_thm_asymp}, as well as for Theorem \ref{nonincr_thm}.
Note that 
$W(\psi,\A) = \limsup_{n \to \infty}A_n$ where
\[A_n := \bigcup_{0 \leq a < a_n} \left[\frac{a}{a_n} - \frac{\psi(n)}{a_n}, \frac{a}{a_n} + \frac{\psi(n)}{a_n}\right].\]
A stronger version of the divergence Borel--Cantelli lemma (see e.g. \cite[Lemma DBC]{BV}) shows that it suffices to establish ``quasi-independence on average'', i.e., 
\begin{equation}\label{QIA}\sum_{m,n \leq N} \lambda(A_n \cap A_m) \ll \sum_{m,n\leq N}\lambda(A_n)\lambda(A_m) \quad\text{ for i.m. } N.\end{equation}
The main obstacle to proving \eqref{QIA} is that a rational number has various representations as a fraction of integers, i.e., it might happen that $\frac{a}{a_n} = \frac{b}{a_m}$ for $m \neq n$.
Thus we write
\[\lambda(A_n \cap A_m) = P(A_n,A_m) + N(A_n,A_m),
\]
where $P(A_n,A_m)$ denotes the contribution of the so-called ``perfect overlaps'' $\frac{a}{a_n} = \frac{b}{a_m}$, while $N(A_n,A_m)$ denotes the contribution of ``non-perfect overlaps'' $\frac{a}{a_n} \neq \frac{b}{a_m}$.
Elementary calculations reveal that indeed $N(A_n,A_m)  \leq \lambda(A_n)\lambda(A_m)$, while $P(A_n,A_m)$ might be too big: 
by a simple Chinese remainder argument,
\[P(A_n,A_m) \asymp \gcd(a_n,a_m)\min\left\{\frac{\psi(n)}{a_n}, \frac{\psi(m)}{a_m}\right\}
\ll \sqrt{\psi(n)}\sqrt{\psi(m)}\frac{\gcd(a_n,a_m)}{\sqrt{a_na_m}}.
\]
Weighted versions of \eqref{def_gn} then imply
\[\sum_{n \neq m\leq N }P(A_n,A_m) \ll \left(\sum_{n \leq N}\lambda(A_n)\right) G(N).\]
Consequently, if $\sum_{n \leq N}\lambda(A_n) \asymp \sum_{n \leq N} \psi(n) \gg G(N)$,
then \eqref{QIA} holds, which leads to \eqref{extra_div_G}.

The estimates on the second moment of the above set system are essentially sharp, so one needs to decrease the actual ``blowup'' of the second moment by working with suitable subsets $A_n' \subseteq A_n$. To avoid 
the perfect overlaps provided by non-uniqueness of fractions, the natural approach is to define

\[A_n' := \bigcup_{\substack{0 \leq a < a_n\\\gcd(a,a_n) =1}} \left[\frac{a}{a_n} - \frac{\psi(n)}{a_n}, \frac{a}{a_n} + \frac{\psi(n)}{a_n}\right].\]
This immediately makes the fractions considered unique, implying for $m \neq n$,

\begin{equation}\label{no_overlap}P(A_n',A_m') = 0.\end{equation}
Note that the estimate $N(A_n',A_m') \ll \lambda(A_n)\lambda(A_m)$
is not always sufficient since in general, 
$\lambda(A_n') \asymp \frac{\varphi(a_n)}{a_n}\lambda(A_n) \not\gg \lambda(A_n)$. Winning back the $\varphi$-factors, i.e., showing on average
\begin{equation}\label{non_perf_DS}N(A_n',A_m') \ll \lambda(A_n')\lambda(A_m'),\end{equation}
 was the substantial technical challenge of the Duffin--Schaeffer Conjecture, which was overcome by Koukoulopoulos--Maynard \cite{KM20}.
However, as mentioned in the introduction, the choice of $A_n'$ is for arbitrary $\A$ not useful. We therefore seek intermediate sets
$A_n'\subseteq A_n^*\subseteq A_n$ so that the perfect overlaps $P(A_n^*,A_m^*)$ are controllable. However, since $A_n^* $ needs to be substantially bigger than $A_n'$, \eqref{no_overlap} is impossible to get, but a reasonable yet controllable contribution might be achievable. Of course, also the non-perfect overlaps $N(A_n^*,A_m^*)$ need to be bounded by the respective $\lambda(A_n^*) \lambda(A_m^*)$, which seems a  plausible task due to \eqref{non_perf_DS}, as long as $A_n'\subseteq A_n^*$ is defined ``multiplicatively''.

Such a set $A_n^{(AEM)}$ in place of $A_n^*$ was found in \cite{AEM}, defined by\footnote{The variable $p$ in this paper always denotes a prime number.} (essentially) 

\begin{equation}\label{eq:sets_AEM}
    A_n^{(AEM)} := \bigcup_{\substack{0 \leq a < a_n\\p \mid \gcd(a,a_n) \implies p > n^2}} \left[\frac{a}{a_n} - \frac{\psi(n)}{a_n}, \frac{a}{a_n} + \frac{\psi(n)}{a_n}\right].
\end{equation}
A quick calculation reveals

\begin{equation}\label{meas_AEM}
\lambda(A_n^{(AEM)}) \geq \lambda(A_n)\prod_{p \leq n^2} \left(1 - \frac{1}{p}\right) \gg \frac{\psi(n)}{\log n}.\end{equation}
Crucially, the cutoff for the sieve condition in \eqref{eq:sets_AEM} depends on $n$ rather than $a_n$, so the
measure bound \eqref{meas_AEM} is uniform in $\A$, as desired.
For the sake of this exposition, we present the ideas for a common cutoff $N^2$ for all $n\leq N$,
suppressing minor inconveniences.
The choice of $N^2$ allows us to treat the perfect overlaps quite directly: 
write $a_n = \ans \cdot \anl$ for the $N^2$-smooth and $N^2$-rough parts of $a_n$ and similarly $a_m = \ams\cdot \aml$. For $n \neq m$, the sieve condition on small primes implies
\[P(A_n^{(AEM)},A_m^{(AEM)}) \neq 0 \implies \ans = \ams, \quad \anl \neq \aml.\]
Since $\anl,\aml$ are both $N^2$-rough, this implies in particular, for $a_n < a_m$, 
\begin{equation}\label{trivial_AEM}
    P(A_n^{(AEM)},A_m^{(AEM)}) \leq \psi(m)\frac{\gcd(a_n,a_m)}{a_m} < \frac{\psi(m)}{N^2},
\end{equation} and summing over $m,n \leq N$ establishes trivially 
\[\sum_{m \neq n \leq N} P(A_n^{(AEM)},A_m^{(AEM)})
\ll \frac{1}{N}\sum_{m \leq N}\psi(m) \ll \sum_{n \leq N}\lambda(A_n^{(AEM)}).
\]
The non-perfect overlaps $N(A_n^{(AEM)},A_m^{(AEM)})$ were treated in \cite{AEM} by an adaptation of the estimates developed for \eqref{non_perf_DS}. We circumvent these technical estimates by appealing to Catlin's conjecture, as discussed later.

Note that in the above, the critically used inequality is \eqref{trivial_AEM}, which is clearly sharp for \textit{fixed} pairs $(a_n,a_m)$; however, we are allowed to average over $1 \leq n,m \leq N$, and can exploit the bilinear structure of
$\{(a_n,a_m): 1 \leq n,m \leq N\}$.
Thus we introduce the parameter $T = T(n)$ and write

\[A_n^* := \bigcup_{\substack{0 \leq a < a_n\\p \mid \gcd(a,a_n) \implies p > T(n)}} \left[\frac{a}{a_n} - \frac{\psi(n)}{a_n}, \frac{a}{a_n} + \frac{\psi(n)}{a_n}\right].\]
Again, assuming for simplicity $T(n) = T(N)$ for all $1 \leq n \leq N$, and writing $a_n = \ans\cdot \anl$ for the $T$-smooth and $T$-rough parts, the analysis from above shows that

\[P(A_n^{*},A_m^{*})
\ll \mathds{1}_{[\ans = \ams]}\varphi{(\ans)}\gcd(\anl,\aml)\min\left\{\frac{\psi(n)}{a_n}, \frac{\psi(m)}{a_m}\right\}.\]
It is this point where we now crucially make use of $(a_n)_n$ being increasing and $\psi$ decreasing: we obtain 
\begin{equation}\label{eq:incr_cruc}\min\left\{\frac{\psi(n)}{a_n}, \frac{\psi(m)}{a_m}\right\} \leq \frac{1}{\sqrt{a_na_m}}\min\{\psi(n),\psi(m)\},\end{equation}
and consequently deduce
\[P(A_n^{*},A_m^{*})
\ll \mathds{1}_{[\ans = \ams]}\frac{\varphi{(\ans)}}{\ans}\frac{\gcd(\anl,\aml)}{\sqrt{\anl\aml}}\min\left\{\psi(n), \psi(m)\right\}.\]
Partitioning the integers into sets $S_d$ with $T$-smooth part $\ans = d$, we obtain

\[\sum_{n \neq m \leq N}P(A_n^*,A_m^*) \leq \sum_{d}\frac{\varphi(d)}{d}\sum_{n,m \in S_d}\frac{\gcd(\anl,\aml)}{\sqrt{\anl\aml}}\min\left\{\psi(n), \psi(m)\right\}. \]
If \begin{equation}\label{gcd_sum_sd}\sum_{n,m \in S_d}\frac{\gcd(\anl,\aml)}{\sqrt{\anl\aml}}\min\left\{\psi(n), \psi(m)\right\} \ll \sum_{n \in S_d} \psi(n)\end{equation}
with the implied constant uniform in $d$, the contribution of perfect overlaps turns out to be negligible.
We note that $(\anl)_{n \in S_d}$ contains $\#S_d \leq N$ distinct $T$-rough numbers. Hence, we choose $T = T(N)$ such that every set $\mathcal{M}$ of at most $N$ distinct $T$-rough numbers satisfies
\begin{equation}\label{rough_gcd_unweighted}\frac{1}{N}\sum_{n,m \in \mathcal{M}}\frac{\gcd(n,m)}{\sqrt{nm}} \ll 1,\end{equation}
which implies for arbitrary non-negative weights $(c_n)_n$ that 
\[\sum_{n,m \in \mathcal{M}}\frac{\gcd(n,m)}{\sqrt{nm}}\min\left\{c_n, c_m\right\} \ll \sum_{n \in \mathcal{M}} c_n.\]
By the above discussion, choosing $T$ such that \eqref{rough_gcd_unweighted} holds implies sufficient control on the contribution from the perfect overlaps.
Analogous to \eqref{meas_AEM}, we obtain
\[\lambda(A_n^*) \gg \frac{\psi(n)}{\log (T(n))},\] and the aim is to beat, even on logarithmic scales, the trivial bound $T = N^2$. This is handled by Theorem \ref{doublelog_thm}, whose strategy is described in the Subsection \ref{subsec:heur_gcd}.

The contribution of non-perfect overlaps $N(A_n^*,A_m^*)$ is circumvented by using Catlin's
conjecture, which is a consequence of the Koukoulopoulos--Maynard theorem
(see \cite[Theorem 2]{KM20}):
writing $\psi_{\A}^+(d) := d\sup_{\substack{n\in\mathbb N\\d\mid a_n}}
     \frac{\psi(n)}{a_n}$, we have
\begin{equation}\label{cat_impl}\sum_{d=1}^{\infty}\frac{\varphi(d)}{d}
     \psi_{\A}^+(d)=\infty \implies  \lambda(W(\psi,\mathcal A))= 1.\end{equation}
To see the idea behind this implication, 
let us write 
\[A_d'(\rho) := \bigcup_{\substack{0 \leq a < d\\(a,d) = 1}} \left[\frac{a}{d} - \rho, \frac{a}{d} + \rho\right],\]
and partition
\[A_n = \mathop{\dot{\bigcup}}_{d \mid a_n} A_{d}'\left(\frac{\psi(n)}{a_n}\right).\]
Grouping the sets indexed by $d$ instead of $n$,
the largest approximation radius around a \textit{reduced} fraction $\frac{a}{d}$
that provides a solution to 
\[\lVert a_n\alpha\rVert \leq \psi(n), \quad d \mid a_n\] is exactly
$\frac{\psi_{\mathcal{A}}^+(d)}{d}$, and hence the Koukoulopoulos--Maynard theorem \cite[Theorem 1]{KM20} applied to 
$\psi_{\A}^+(d)$ shows \eqref{cat_impl}.

We transfer this idea from 
$W(\psi,\A) = \limsup_{n \to \infty} A_n$ to
 $\limsup_{n \to \infty} A_n^*$: we partition
\begin{equation*}A^{*}_n = \mathop{\dot{\bigcup}}_{d\mid^* a_n} A_{d}'\left(\frac{\psi(n)}{a_n}\right), \end{equation*}
where (recall $a_n = \ans \anl$ for the $T$-smooth and $T$-rough parts)
\[d \mid^* a_n :\Leftrightarrow d = e\ans\quad \text{ for some }\quad e \mid \anl.\] 
Therefore, following the above discussion, the Koukoulopoulos--Maynard Theorem implies
\[\sum_{d=1}^{\infty}\frac{\varphi(d)}{d}
     \psi_{\A}^*(d)=\infty \implies \lambda(\limsup_{n \to \infty}A_n^*) = 1,\]
where
\begin{equation}\label{def_*}\psi_{\mathcal{A}}^*(d) := d\sup_{\substack{n: d\mid^* a_n}}
     \frac{\psi(n)}{a_n}.\end{equation}
Thus, one is left to show that
\[\sum_{d}\frac{\varphi(d)}{d}
     \psi_{\A}^*(d) \asymp \sum_n \lambda(A_n^*) \geq \sum_n \frac{\psi(n)}{\log T(n)} \to \infty.\]
This reduces to showing that, in an averaged sense,
\[\sup_{\substack{n: d\mid^* a_n}}
     \frac{\psi(n)}{a_n} \gg \sum_{n: d\mid^* a_n}\frac{\psi(n)}{a_n},\]
or in other words,
\begin{equation}\label{avg_sup_sum}\sum_d \varphi(d)\sum_{n\in S_d}\frac{\psi(n)}{a_n}
\ll \sum_d \frac{\varphi(d)}{d}\psi^*_{\A}(d).
\end{equation}
Clearly, for $f \geq 0$,
\[\left(\sum_{x \in X}f(x)\right)^2 \ll (\sup_{x \in X}f(x))\cdot  \sum_{x \in X}f(x) \implies
\sum_{x \in X}f(x) \ll \sup_{x \in X}f(x),
\]
so in order to show \eqref{avg_sup_sum}, we use Cauchy--Schwarz which gives
\[
\begin{split}
\left(\sum_d \varphi(d)\sum_{n\in S_d}\frac{\psi(n)}{a_n}\right)^2
&\leq
\left(\sum_d \frac{\varphi(d)}{d}\psi_{\A}^*(d)\right)
\left(\sum_d \varphi(d)
      \sum_{n,m\in S_d}\min\left\{\frac{\psi(n)}{a_n},\frac{\psi(m)}{a_m}\right\}\right).
\end{split}
\]
The second factor above equals \textit{exactly} the contribution of the perfect overlaps
$\sum_{n ,m}P(A_n^*, A_m^*).$
This is not a coincidence, but rather arises from the construction of $\psi^*_{\A}$. From the previous discussion, the terms
$\sum_{n ,m}P(A_n^*, A_m^*)$ are controllable once the respective GCD sums are understood, and the Catlin criterion allows us to bypass controlling the non-perfect overlap contribution completely.\footnote{Handling non-perfect overlaps is the crucial technical achievement by Koukoulopoulos--Maynard. Hence the proof of our results still relies on it, but a lengthy direct computation is circumvented by applying the ``blackbox theorem'' of Catlin's conjecture.}

Let us finally reflect on why this approach goes (provably) wrong for arbitrary integer sequences. We used the property of $(a_n)_n$ being increasing in \eqref{eq:incr_cruc};
if we worked instead with the always valid
\begin{equation}\label{agm}\min\left\{\frac{\psi(n)}{a_n}, \frac{\psi(m)}{a_m}\right\}
\leq \frac{\sqrt{\psi(n)\psi(m)}}{\sqrt{a_na_m}},\end{equation}
one would need to show that any $T$-rough set $\mathcal{M}$ satisfies
\begin{equation}\label{spectral}\sum_{n,m \in \mathcal{M}}\frac{\gcd(n,m)}{\sqrt{nm}}c_n^{1/2} c_m^{1/2} \ll \sum_{n \in \mathcal{M}} c_n\end{equation}
instead of
\begin{equation}\label{min_form}\sum_{n,m \in \mathcal{M}}\frac{\gcd(n,m)}{\sqrt{nm}}\min\left\{c_n, c_m\right\} \ll \sum_{n \in \mathcal{M}} c_n,\end{equation}
as used above. If we could find $T = N^{o(1)}$ such that \eqref{spectral} holds for any $T$-rough set $\mathcal{M}$ and arbitrary weights $(c_n)_n > 0$, then this would lead to an improvement over the rate $\frac{1}{n}$ for arbitrary sequences $\A$ in Theorem \ref{nonincr_thm},
which by \eqref{eq:1n_sharp} is ruled out. Hence \eqref{spectral} cannot hold for any $T = N^{o(1)}$, in contrast to the $T$ for \eqref{min_form} provided in Theorem \ref{doublelog_thm}.
The example we provide for \eqref{eq:1n_sharp} is built on the Duffin--Schaeffer counterexample \cite{DS41}, and indeed fails the above estimate \eqref{spectral} for any $T = N^{o(1)}$. As in \cite{DS41}, we consider large products of primes, and include all its divisors in the support of $(a_n)_n$. By ordering the sequence accordingly, we make sure that $\frac{\psi(n)}{a_n} =\frac{\psi(m)}{a_m}$ for many distinct $m,n$ (compare \eqref{agm}). This particular divisor structure as well as the particularly ``maliciously chosen'' weights creates too much overlaps for the $\limsup$-set to have positive measure.

\subsection{Ideas for the uniform approximation}\label{subsec:idea_uni}

A direct treatment of the questions concerning \textit{uniform} approximation with second moment methods described above is very inefficient: the approximation functions $\psi_N(n), 1 \leq n \leq N$ depend on $N$, and thus a na\"{\i}ve second moment and Chebyshev argument would need to sum the exceptional probabilities over \textit{all} $N$. We reduce the problem to studying \textit{quantitative asymptotic} approximations: we show that ``with high probability'', a fixed dyadic range $[2^j,2^{j+1})$ contains at least one 
$n \in [2^j,2^{j+1})$ that satisfies $\lVert a_n\alpha\rVert \leq \psi(n)$. By a sandwiching argument that only loses a constant factor, this implies the uniform  statement. 
To obtain this quantitative control, it suffices to obtain good bounds on the variance of the counting function
\[\mathcal{N}_j(\alpha) := \# \{n \in [2^j,2^{j+1}): \lVert a_n \alpha\rVert \leq\psi(n)\},\]
which is expected to concentrate around $\Psi_j := \sum_{n \in [2^j,2^{j+1})} 2\psi(n)$.
In this way, we have reduced the uniform problem to an asymptotic quantitative one, sometimes known as ``Schmidt-type'' due to Schmidt's influential article \cite{Schmidt} in the Khintchine setup (with monotonicity). In the Duffin--Schaeffer setup, i.e., $\mathcal{N}_j'$ only counting approximations 
\[\lvert a_n\alpha - a| \leq \psi(n), \quad \gcd(a,a_n) = 1,\, n \in [2^j,2^{j+1}),\] quantitative theory for arbitrary functions was recently developed \cite{ABH,HSW24,KMY24}. The sharpest established version (\cite[Theorem 2]{KMY24}) comes with an almost square-root cancellation in the arising variance calculation. This would allow to deduce by Chebyshev's inequality,
\[\lambda(\{\alpha: \mathcal{N}_j(\alpha) = 0 \}) \leq \lambda(\{\alpha: |\mathcal{N}_j'(\alpha) - \Psi_j'| \geq \Psi_j'/2 \}) \ll_{\varepsilon} \frac{1}{(\Psi_j')^{1 - \varepsilon}},
\]
where $\Psi_j' = \sum_{n \in [2^j,2^{j+1})} 2\frac{\varphi(a_n)}{a_n}\psi(n)$.
This permits a Borel--Cantelli argument once $\Psi_j' \gg j^{1 + 2\varepsilon}$,
which exactly shows the reason we need to have this $(\log N)^{1 + \varepsilon}$ term to ``spend''.

Of course, again, fully-reduced fractions would lose the unbounded factor $\frac{\varphi(a_n)}{a_n}$, hence we work again with partially reduced sets $A_n^*(N), 1 \leq n \leq N$, 
and define $\psi_{\A,N}^*(d)$ similarly to \eqref{def_*}. By a Catlin-type partition as described in Section \ref{subsec:asymp}, again using the $T$-roughness property to take care of perfect overlaps,
 we are left to count solutions to
\[\lvert d\alpha - a\rvert \leq \psi_{2^j}^*(d), \quad \gcd(a,d) = 1,\]
which can be handled by the quantitative Duffin--Schaeffer theory described above.

We remark here that the best available variance estimate \cite[Theorem 2]{KMY24} comes with exponent $\frac{1}{2} + \varepsilon$, instead of a potential $\frac{1}{2} + o(1)$. Thus, the result in Theorem \ref{main_thm_unif} 
is unaffected by the improvement from our more involved GCD sum approach over taking the trivial roughness parameter $T = N^2$: the improvement obtained by Theorem \ref{doublelog_thm} is (after taking $\log$) on a $\log \log$-scale, and Theorem \ref{unbounded_GCD} shows that every potential improvement is at most
$(\log N)^{o(1)}$; this gain is absorbed in the $(\log N)^{\varepsilon}$-term that already needs to be ``spent'' for the $\frac{1}{2} + \varepsilon$ in \cite[Theorem 2]{KMY24}. Hence, the methods developed in this paper do not allow to improve the current rate of Theorem \ref{main_thm_unif} when assuming $\A$ to be increasing, in contrast to the improvement gained in the asymptotic approximation (cf. Theorem \ref{nonincr_thm} vs. Theorem \ref{main_thm_asymp}).

\subsection{Ideas for results about GCD sums}\label{subsec:heur_gcd}

For Theorem \ref{doublelog_thm}, we partially follow the strategy employed by Bondarenko and Seip \cite{BS_15}: we first establish an upper bound for more general \textit{weighted} GCD sums over sets of \textit{square-free} integers.
 An idea that goes back to G\'{a}l \cite{Gal} shows that one may restrict attention to divisor-closed sets. One key argument of Bondarenko--Seip is to further show that maximizing GCD sums arise for so-called complete sets. This completeness property imposes strong combinatorial restrictions on the occurrence of large prime factors, since their presence forces many additional integers to belong to the set. Taking into account that all prime divisors of $T$-rough sets are larger than $T$, we introduce an analogue of this property, namely the T-completeness property, which imposes restrictions on the number of (even) large(r) prime factors. Then, the proof of Theorem \ref{doublelog_thm} deviates from earlier works: we partition the pairs $(m,n) \in \mathcal{M}^2$ into two subsets depending on whether or not $m$ and $n$ share all their prime factors up to $V$ (where $T<V = T+ o(T)$). If they do, the divisor-closed and completeness property implies that these pairs differ only on a small set of very large prime factors, and their contribution can be bounded using a sieve argument. If there exists $p\in (T, V)$ such that $p$ divides exactly one of $n,m$, we bound their contribution using the standard GCD sum bound \eqref{gcd_asymp}.  \\

For Theorem \ref{unbounded_GCD}, we construct sets $\mathcal{M}$ consisting of $T$-rough integers with large associated GCD sums. To do so, we choose sets $\mathcal{M}$ of $N$ squarefree numbers with a fixed number $k \asymp (\log N)^{\varepsilon}$ of prime divisors in an interval $(T,U)$ with $\frac{\log U}{\log T} \asymp (\log N)^{2\varepsilon}$. This construction is suggested by the approach of de la Bretèche--Tenenbaum \cite{BT}. However, directly translating the ensemble from \cite{BT} to our setup would only allow us to construct large GCD sums for $T= \exp((\log N)^{1/2-\varepsilon})$. In order to obtain our result with $T=\exp((\log N)^{1-\varepsilon})$, we impose a minimality condition by defining $\mathcal{M}$ as the $N$ smallest integers satisfying these properties. These numbers, which are all $\leq N^{1+\varepsilon}$, have unusually many prime factors compared to a generic integer of this size. 
Given $I$ a subset of prime divisors of $n$, and $p_I$ the product of those primes, we only consider the contribution of these $m \in \mathcal{M}$ where $m=\frac{n}{p_I} q$ with $q \leq p_I$. The contribution of such a pair $(m,n)$ then gives $\frac{\gcd(m,n)}{\sqrt{mn}}= \frac{1}{\sqrt{p_I}}\frac{1}{\sqrt{q}}$. 
For $j>1/\varepsilon$, we establish that for most $n \in \mathcal{M}$, we can construct many different subsets $\vert I\vert =j $ such that $p_I$ is not too large. Averaging over these subsets, we prove (using the minimality condition) that these pairs give a large contribution to the GCD sum.

\section{Proofs for bounded and unbounded GCD sums}\label{sec:GCD}

\subsection{The proof of Theorem \ref{doublelog_thm} about bounded GCD sums}

\subsubsection*{Prerequisites}

In this section we gather various estimates that are known from earlier works on upper bounds for GCD sums. While classical upper bounds on GCD sums do not consider a roughness assumption, we still make use of several auxiliary results developed there, and adapt them to the $T$-rough setting. Note that if $f$ denotes the multiplicative function such that
$f(p) = \frac{1}{\sqrt{p}}$, we can write 

\[\sum_{n,m \in \mathcal{M}} \frac{\gcd(m,n)}{\sqrt{mn}}
= \sum_{n,m \in \mathcal{M}} f\left(\frac{mn}{\gcd(m,n)^2}\right) =: \mathcal{S}(\mathcal{M},f). 
\]
For technical reasons, we need to work with a family of more general completely multiplicative functions, defined by $g(p) := \min\{\frac{C}{\sqrt{p}},1\}$. Note here that (for $T > C$)
\[\mathcal{S}(\mathcal{M},g) := \sum_{m,n \in \mathcal{M}} g\left(\frac{mn}{\gcd(m,n)^2}\right) =
\mathcal{S}(\varphi_T(\mathcal{M}),g_T),\]
where $\varphi_T: \{N \in \mathbb{N}: p \mid N \implies p > T\} \to \N$ is defined as the completely multiplicative function with\footnote{as usual, $p_n$ denotes the $n$-th smallest prime, and $\pi(x) = |\{p \leq x: p \in \mathbb{P}\}|$ is the prime counting function.}
$\varphi_T(p_n) = p_{n-\pi(T)}$  and $g_T(p_n) := \frac{C}{\sqrt{p_{n+\pi(T)}}}$. Since $\varphi_T(\mathcal{M})$ is a set of cardinality $N$, and $(g_T(p))_{p \in \mathbb{P}}$ is a decreasing sequence, we obtain from \cite[Lemmas 1 - 3]{BS_15} and \cite[Lemma 2]{ABS_15} the statements listed below.\\

A set $\mathcal{M}$ is called \textit{divisor-closed} whenever $n \in \mathcal{M} \implies d \in \mathcal{M}$ for all $d \mid n$. We use the following variant of \cite[Lemma 1]{BS_15}.

\begin{lem}\label{lem1_bs}
Let $g$ be completely multiplicative, with $g(p) := \min\{\frac{C}{\sqrt{p}},1\}$. Then there exists a divisor-closed $T$-rough set $\mathcal{M}$ consisting of $N$ distinct square-free integers such that
  \begin{equation}\label{extremal_sf}\mathcal{S}(\mathcal{M},g) = \sup_{\substack{\mathcal{M}'\, T-\text{rough}\\|\mathcal{M'}| = N\\\mathcal{M}' \text{square-free}}}\mathcal{S}(\mathcal{M}',g).\end{equation}
\end{lem} Let $\tau(n)$ denote the number of divisors of $n$, and $\omega(n)$ the number of distinct prime divisors of $n$.
We note for further reference that \begin{equation}\label{div_bound_triv}\mathcal{M} \text{ divisor-closed} \implies \forall n \in \mathcal{M}, \quad \omega(n) \leq  \frac{\log(\lvert \mathcal{M}\rvert)}{\log 2},\end{equation}
which follows immediately from $2^{\omega(n)} \leq \tau(n) \leq \vert \mathcal{M}\vert $.\\

In analogy to the \textit{completeness property} in \cite{BS_15}, we define the following \textit{$T$-completeness property}: a set $\mathcal{M}$ of distinct positive $T$-rough integers has the $T$-completeness property if for every $n \in \mathcal{M}$ and prime $p$ such that $p\mid n$, we have for all primes $T < p' < p$ that $p'\mid n$ or $n p^{-1}p' \in \mathcal{M}$.
This definition allows us to obtain the following variant of \cite[Lemma 2]{BS_15}.

\begin{lem}\label{lem2_bs}
Let $\mathcal{M}$ be a divisor-closed set of square-free numbers satisfying \eqref{extremal_sf}. Then 
$\mathcal{M}$ is $T$-complete.
\end{lem}
We rephrase \cite[Lemma 2]{ABS_15} (see also \cite[Lemma 5]{BS_15}).
\begin{lem}\label{lem2_abs}
Let $f,g$ be completely multiplicative functions, defined by $f(p) = \frac{1}{\sqrt{p}}, g(p) = \min\{\frac{2}{\sqrt{p}},1\}$.
Then for every $T$-rough set $\mathcal{M}$ of $N$ distinct integers, there exists a $T$-rough divisor-closed set $\mathcal{M}'$ consisting of $N$ distinct integers such that
\[\mathcal{S}(\mathcal{M},f) \leq \mathcal{S}(\mathcal{M}', g).\]
\end{lem}

Next, we adapt \cite[Lemma $3$]{BS_15} (see also \cite[Lemme $3.1$]{BT}) to the $T$-rough setup considered here.
\begin{lem}\label{lem3_bs}
    Let $\mathcal{M}$ be a $T$-complete and divisor-closed set of $N$ distinct integers. Let $V>T$ such that $U := |\{T < p < V\}| \geq 4(\log N)^2$.  Then for all $m \in \mathcal{M}$, we have
\[|\{p\geq V: p \mid m\}| \leq 4\frac{\log N}{\log U}.\]
\end{lem}

\begin{proof}
Let $p_{j_1}<\dots<p_{j_r}$ be the prime factors of $m$ exceeding $V$. By the $T$-completeness and divisor-closed properties of $\mathcal{M}$, there exist at least $\left(\prod_{k=1}^{r}(U-k+1)\right)/r!$ integers in $\mathcal{M}$ of the form $\prod_{1 \leq k \leq r} p_{s_k}$ with $\pi(T)<s_k\leq j_r.$ By \eqref{div_bound_triv} and our assumption on $U$ we have $r! \leq \left(\frac{\log N}{\log 2}\right)^r \leq U^{r/2}$. Moreover, $$U-(k-1)= U\left(1-\frac{(k-1)}{U}\right) \geq U\left(1-\frac{1}{\sqrt{U}}\right), \qquad k=1,\dots,r.$$ Thus, from $r \leq \sqrt{U}$,
$$ N \geq \frac{\left(\prod_{k=1}^{r}(U-k+1)\right)}{r!} \geq U^{r/2} \left(1-\frac{1}{\sqrt{U}}\right)^r \geq \frac{U^{r/2}}{4}.$$ Taking logarithms, the conclusion follows.
\end{proof}

\begin{lem}\label{lem:weighted_GCD}
Let $g$ be a non-negative function, and assume $\mathcal{M} \subseteq \N$ to be such that for all nonempty $\mathcal{M}' \subseteq \mathcal{M}$, 
\[\mathcal{S}(\mathcal{M}',g) \ll |\mathcal{M}'|,\]
with a uniform implied constant. 
Then for any non-negative weights $(c_n)_{n \in \mathcal{M}}$, we have
    \[\sum_{n,m\in \mathcal{M}} g\left(\frac{mn}{\gcd(m,n)^2}\right)\min\{c_n,c_m\} \ll \sum_{m \in \mathcal{M}}c_m.\]
\end{lem}
\begin{proof}
Without loss of generality, we may assume that $c_m \leq 1$ for all $m$, otherwise we renormalize. 
    Given $i \in \mathbb{N}$, we write
    $\mathcal{M}_i := \{m \in \mathcal{M}: c_m \geq 2^{-i}\}$.
    We then obtain
    \[\sum_{n,m\in \mathcal{M}} g\left(\frac{mn}{\gcd(m,n)^2}\right)\min\{c_n,c_m\}
    \ll \sum_{i \geq 0} 2^{-i}\sum_{n,m\in \mathcal{M}_i} g\left(\frac{mn}{\gcd(m,n)^2}\right).
    \]
    Since $\mathcal{M}_i \subseteq \mathcal{M}$, we get
    $\sum_{n,m\in \mathcal{M}_i} g\left(\frac{mn}{\gcd(m,n)^2}\right)
    \ll \vert\mathcal{M}_i\vert
    $ and consequently,
    \[\begin{split}\sum_{n,m\in \mathcal{M}} g\left(\frac{mn}{\gcd(m,n)^2}\right)\min\{c_n,c_m\} &\ll \sum_{i \geq 0} 2^{-i} \vert \mathcal{M}_i\rvert
    \ll \sum_{i \geq 0}2^{-i} \sum_{j \leq i} |\{m \in \mathcal{M}: c_m \asymp 2^{-j}\}|
    \\&\ll \sum_{j \geq 0} |\{m \in \mathcal{M}: c_m \asymp 2^{-j} \}|\sum_{i \geq j}2^{-i}\ll \sum_{m \in \mathcal{M}}c_m.
    \end{split}\]
\end{proof}

The following contains the technical novel part of the proof of Theorem \ref{doublelog_thm}, having adapted the necessary combinatorial and number-theoretic properties above.
\begin{lem}\label{lem:square_free_gen}
    Let $\mathcal{M}$ be a set of at most $N$ squarefree $T$-rough numbers.
    Let $g$ denote the completely multiplicative function defined by
    $g(p) = \min\{2/\sqrt{p},1\}$.
    Then 
    \[\mathcal{S}(\mathcal{M},g) \ll |\mathcal{M}|.\]
\end{lem}

\begin{proof}
    By Lemmas \ref{lem1_bs} and \ref{lem2_bs}, we may assume without loss of generality that $\mathcal{M}$ is divisor-closed and $T$-complete. 
    Let $V = T\left(1 + \tfrac{1}{\exp((\log N)^{4/7})}\right)$ and set
 \[\mathcal{A} := \{(n,m) \in \mathcal{M} \times \mathcal{M}: \exists p \in (T, V): \nu_p(n) \neq \nu_p(m)\},\]
    and 
    \[\mathcal{B} := \{(n,m) \in \mathcal{M} \times \mathcal{M}: \forall p \in (T, V): \nu_p(n) = \nu_p(m)\},\] where $\nu_p$ denotes the $p$-adic valuation.
    With $\varphi_T$ being the prime-shift defined in the beginning of this section, we write
\[\sum_{(m,n) \in \mathcal{A}} g\left(\frac{mn}{\gcd(m,n)^2}\right) = \sum_{(m,n) \in \mathcal{A}} g\left(\frac{\varphi_T(m)\varphi_T(n)}{\gcd(\varphi_T(m),\varphi_T(n))^2}\right)
\frac{g\left(\frac{mn}{\gcd(m,n)^2}\right)}{g\left(\frac{\varphi_T(m)\varphi_T(n)}{\gcd(\varphi_T(m),\varphi_T(n))^2}\right)}.\]
If $(m,n) \in \mathcal{A}$, by definition, there must be at least one prime $T < \hat{p} = \hat{p}(m,n) < V$ where $\nu_{\hat{p}}(m) \neq \nu_{\hat{p}}(n)$. Since $g(p)_{p \in \mathbb{P}}$ is monotonically decreasing, we obtain

\begin{align*}\frac{g\left(\frac{mn}{\gcd(m,n)^2}\right)}{g\left(\frac{\varphi_T(m)\varphi_T(n)}{\gcd(\varphi_T(m),\varphi_T(n))^2}\right)}&
\leq \frac{g\left(\hat{p}\right)}{g(\varphi_T(\hat{p}))} = \sqrt{\frac{\varphi_T(\hat{p})}{\hat{p}}}  
\ll \sqrt{\frac{(V-T) \log(V-T)}{T}} \leq \frac{1}{\exp(\tfrac{1}{3}(\log N)^{4/7})}.
\end{align*}
Here we used $\pi(\hat{p})-\pi(T)\leq V-T+1$ and
$p_k\ll k\log k$. We get
\[\begin{split}\sum_{(m,n) \in \mathcal{A}} g\left(\frac{mn}{\gcd(m,n)^2}\right)
&\leq \frac{1}{\exp(\tfrac{1}{3}(\log N)^{4/7})}\sum_{(m,n) \in \mathcal{A}} g\left(\frac{\varphi_T(m)\varphi_T(n)}{\gcd(\varphi_T(m),\varphi_T(n))^2}\right)
\\&\ll \frac{1}{\exp(\tfrac{1}{3}(\log N)^{4/7})}\sum_{m,n \in \varphi_T(\mathcal{M})} g\left(\frac{mn}{\gcd(m,n)^2}\right)
\\&= \frac{1}{\exp(\tfrac{1}{3}(\log N)^{4/7})}\mathcal{S}(\varphi_T(\mathcal{M}),g)
\\&\ll \frac{|\mathcal{M}| \exp((\log N)^{1/2 + 1/100})}{\exp(\tfrac{1}{3}(\log N)^{4/7})} \ll |\mathcal{M}|,
\end{split}
\]
where we used \cite[Theorem 2]{BS_15} (which is a more general version of the GCD sum bound \eqref{gcd_asymp})
to guarantee $\mathcal{S}(\varphi_T(\mathcal{M}),g) \ll |\mathcal{M}|\exp((\log N)^{1/2 + 1/100})$.

We now bound the contribution of $\mathcal{B}$. For $n \in \mathcal{M}$, let $$n_s = \prod_{T < p < V}p^{\nu_p(n)} \qquad \text{and } \qquad n_r = \prod_{p \geq V}p^{\nu_p(n)}$$ denote the $V$-smooth and $V$-rough part of $n$.
We then define
\[\mathcal{B}_n := \{m \in \mathcal{M}: (m,n) \in \mathcal{B}\} = \{m \in \mathcal{M}: m_s = n_s\}.\]
Hence we obtain

\begin{equation}\label{contr_B}\sum_{(n,m) \in \mathcal{B}}g\left(\frac{nm}{\gcd(n,m)^2}\right)
\leq |\mathcal{M}| + 2 \sum_{n \in \mathcal{M}}
\sum_{\substack{m \in \mathcal{B}_n\\ m < n}}
g\left(\frac{n_rm_r}{\gcd(n_r,m_r)^2}\right).
\end{equation}
By the prime number theorem version of Vinogradov--Korobov 
\[\pi(x) = \operatorname{Li}(x) + O(x\exp(-c(\log x)^{3/5-1/100})),\] we get that 
\[U := \#\{T < p < V\} \gg \frac{(V-T)}{\log T}.\]
Since $\mathcal{M}$ is divisor-closed and $T$-complete, Lemma \ref{lem3_bs} proves that for any $m \in \mathcal{B}_n$,
\begin{equation}\label{bnd_primediv}\omega(m_r) \leq \frac{4\log N}{\log U} \leq \frac{5\log N}{\log T} = \frac{\log \log N}{4}.\end{equation} 
Therefore, we get
\begin{equation}\label{g_to_sqrt}g\left(\frac{n_rm_r}{\gcd(n_r,m_r)^2}\right) \leq 2^{\omega(n_rm_r)}\frac{\gcd(n_r,m_r)}{\sqrt{n_rm_r}} \leq \frac{\gcd(n_r,m_r)}{\sqrt{n_rm_r}}(\log N)^{(\log 2)/2}.\end{equation}
Note that \[\begin{split}
 \frac{1}{\sqrt{n_r}} \sum_{\substack{m \in \mathcal{B}_n\\ m < n}}
\frac{\gcd(n_r,m_r)}{\sqrt{m_r}}  
& \leq \frac{1}{\sqrt{n_r}}\sum_{\substack{d \mid n_r\\d \neq n_r}}\sqrt{d}\sum_{\substack{j \leq n_r/d\\j\, V-\text{rough}
}} \frac{1}{\sqrt{j}},
\end{split}
\]
where we have parametrized $d = \gcd(m_r,n_r)$ and $m_r = j \gcd(m_r,n_r)$.
By a straightforward application of an upper-bound sieve (see e.g. \cite[Theorem 18.11]{kouk_book}),
we get (note that $n_r/d \geq V$)
\[\sum_{\substack{j \leq n_r/d\\j\, V-\text{rough}}}\frac{1}{\sqrt{j}}
 \ll \frac{\sqrt{n_r/d}}{\log V} \leq   \frac{\sqrt{n_r/d}}{\log T}.
\]
Thus, we obtain 

\[\frac{1}{\sqrt{n_r}} \sum_{\substack{m \in \mathcal{B}_n\\ m < n}}
\frac{\gcd(n_r,m_r)}{\sqrt{m_r}}  
\ll \frac{1}{\log T}\sum_{d \mid n_r}1
\ll \frac{ (\log N)^{(\log 2)/4}}{(\log N)^{1-o(1)}},\] since (recall that $n_r$ is squarefree) $\tau(n_r) = 2^{\omega(n_r)}  \leq (\log N)^{(\log 2)/4}$ by \eqref{bnd_primediv}. Summing over $n \in \mathcal{M}$, using \eqref{contr_B} and \eqref{g_to_sqrt} and combining with the contribution of $\mathcal{A}$, we obtain

\[\sum_{m,n \in \mathcal{M}}g\left(\frac{mn}{\gcd(m,n)^2}\right)
\ll \lvert \mathcal{M}\rvert,
\]
which finishes the proof.
\end{proof}

We finally connect all pieces from above to the complete proof.
\begin{proof}[Proof of Theorem \ref{doublelog_thm}]
   By Lemma \ref{lem2_abs} (using $f,g$ as in the lemma),
\[\sum_{m,n \in \mathcal{M}}\frac{\gcd(m,n)}{\sqrt{mn}} = 
\sum_{m,n \in \mathcal{M}}f\left(\frac{nm}{\gcd(m,n)^2}\right)
\ll 
\sum_{m',n'\in \mathcal{M}'}g\left(\frac{n'm'}{\gcd(m',n')^2}\right),
\]
where $\mathcal{M}'$ is a divisor-closed $T$-rough set with $|\mathcal{M'}| = N$, and $g(p) :=  \min\{2/\sqrt{p},1\}$.
Let $\mathcal{M}'' = \rad(\mathcal{M}')$ where $\rad(n) := \prod_{p \mid n} p$ denotes the radical of $n$. For $n \in \mathcal{M}''$, we define 
\[\mathcal{S}_n := \{m' \in \mathcal{M}': \rad(m') = n\}.\]
Let $n,m \in \mathcal{M}''$. Then for $|\mathcal{S}_n| \leq |\mathcal{S}_m|$, we have by following the argument in \cite[proof of Theorem 1., p.38]{BS_15} that

\[\begin{split}\sum_{n' \in \mathcal{S}_n}\sum_{\substack{m' \in \mathcal{S}_m}}
g\left(\frac{{n'm'}}{\gcd(n',m')^2}\right) 
&\ll g\left(\frac{n m }{\gcd(n,m)^2}\right)
\sum_{n' \in \mathcal{S}_n}
\prod_{p\mid m}\left(1 - \frac{4}{\sqrt{p}}\right)^{-1}.
\end{split}
\]
Note that $\mathcal{M}'$ is divisor-closed and hence by \eqref{div_bound_triv} satisfies $\omega(n') \leq \frac{\log N}{\log 2}, n' \in \mathcal{M}'$. Thus by $T$-roughness, we have
\[\prod_{p\mid m}\left(1 - \frac{4}{\sqrt{p}}\right)^{-1} \ll \exp\left(\sum_{T<n\leq T+ \log N/\log 2} \frac{4}{\sqrt{n}}\right) \ll 1,\]
showing that 
\begin{align*}\label{boundSnSm} \sum_{m,n \in \mathcal{M}}\frac{\gcd(m,n)}{\sqrt{mn}}& \ll \sum_{n,m \in \mathcal{M}''} \left(\sum_{n' \in \mathcal{S}_n}\sum_{m' \in \mathcal{S}_m}
 g\left(\frac{{n'm'}}{\gcd(n',m')^2}\right) \right) \\
& \ll  \sum_{n,m \in \mathcal{M}''} g\left(\frac{n m}{\gcd(n,m)^2}\right) \min\{|\mathcal{S}_n|,|\mathcal{S}_m|\}.\end{align*}
In view of Lemma \ref{lem:weighted_GCD}, it therefore suffices to show that
for every $\mathcal{M}'''\subseteq \mathcal{M}''$, we have
\[\mathcal{S}(\mathcal{M}''',g) \ll |\mathcal{M}'''|.\]
Note that $\mathcal{M}''$ is squarefree, thus so is $\mathcal{M}'''$. Hence the unweighted statement \eqref{eq:unw_bd} follows from Lemma \ref{lem:square_free_gen}.
The weighted statement \eqref{eq:weight_bd}, and thus Theorem \ref{doublelog_thm}, follows from yet another application of Lemma \ref{lem:weighted_GCD}.
\end{proof}

\subsection{The proof of Theorem \ref{unbounded_GCD} about unbounded GCD sums}

    By relabeling, it suffices to show Theorem \ref{unbounded_GCD} for $T = \exp(\log(N)^{1-3\varepsilon})$, and by monotonicity of the roughness condition, we can further assume that $\varepsilon = \frac{1}{A}$ for some positive integer $A \geq 4$.
We assume $N$ to be sufficiently large and let \begin{equation}\label{defU}U = \exp((\log N)^{1 - \varepsilon}), \qquad  k = \lfloor(\log N)^{2\varepsilon}\rfloor. \end{equation} 
We now define our set as\footnote{$\mu$ denotes the Möbius function.}\footnote{The interested reader might notice (this can be easily computed following the lines of \cite[Section $4$]{Mun-energy}) that for $\log T \geq (\log N)^{1/2+\varepsilon}$, the set of the first $T$-rough numbers has its associated GCD sums bounded.} 
\begin{equation}\label{defM}\mathcal{M} := \{\text{the $N$ smallest numbers $m$: } \omega(m) = k, \mu^2(m) = 1, p\mid m \implies p \in [T,U]\}.\end{equation}
We see that this set is well-defined by the fact that $\binom{|\mathbb{P} \cap [T,U]|}{k} \gg \frac{U^k}{(k\log U)^k} \geq N$ for $N$ sufficiently large by the prime number theorem.\\

For $n \in \mathcal{M}$, let $\mathcal{P}(n) = \{p \in \mathbb{P}: p\mid n\}$ be the set of prime divisors of $n$. For $I \subseteq \mathcal{P}(n)$, 
let $p_I = p_I(n) := \prod_{p \in I}p$. We now define
 \[ \mathcal{M}_{n,I}= \{q \in \mathbb{N},q \leq p_I, (q,n)=1, \omega(q)=|I|, \mu^2(q)=1, p \mid q \implies p \in [T,U] \}.\] For $I \subseteq \mathcal{P}(n)$ and $q \in \mathcal{M}_{n,I}$, we further let \[m_{I,q} = m_{I,q}(n) :=  \frac{n}{p_I} q.  \] Note that $m_{I,q} \in \mathcal{M}$ when $q \in \mathcal{M}_{n,I}$.
Let $j \in \N$, and define the set of $j$-subsets of $\mathcal{P}(n)$ as $$[ \mathcal{P}(n) ]^{j} := \{ I \mid I \subseteq \mathcal{P}(n) \text{ and } | I | = j  \}.$$
For $n = p_1 \cdots p_k$ and $1 \leq j \leq k$, we let 
\[ \mathcal{Q}_{n,j} = \left\{ (I, q) \;\middle\vert{}\; I \in [ \mathcal{P}(n) ]^{j},\; q \in \mathcal{M}_{n,I} \right\}.\] 
The coprimality condition ensures that the following map is injective. \begin{equation}\label{map_phi}
\begin{array}{rcl}
\phi_{n,j} : \mathcal{Q}_{n,j}  & \longrightarrow & \mathcal{M} \\
(I,q) & \longmapsto & m_{I,q} = \frac{n}{p_I}q.
\end{array}
\end{equation} Hence, for a given $n\in \mathcal{M}$ and $j\leq k$ we have \begin{equation}\label{pw_GCD}\sum_{m \in \mathcal{M}} \frac{\gcd(m,n)}{\sqrt{mn}} \geq \sum_{m \in \phi_{n,j}(\mathcal{Q}_{n,j})} \frac{\gcd(m,n)}{\sqrt{mn}}.\end{equation}
Note that the pair $(n,m_{I,q})$ gives a contribution to the GCD sum of
\[\frac{\gcd(m_{I,q},n)}{\sqrt{m_{I,q}n}}
= \frac{1}{\sqrt{p_I}}\frac{1}{\sqrt{q}}.\]
Thus, for every fixed $n \in \mathcal{M}$ and any fixed $j \leq k$, we obtain using \eqref{pw_GCD} that 

\begin{equation}\label{pw_GCD_bis}\sum_{m \in \mathcal{M}} \frac{\gcd(m,n)}{\sqrt{mn}} \geq \sum_{I \in [\mathcal{P}(n)]^j} \frac{1}{\sqrt{p_I}}\sum_{q \in \mathcal{M}_{n,I}} \frac{1}{\sqrt{q}}.
\end{equation}

We now specify an integer $j \in \N$ (large, but independent of $N$), and assume that $j$ is divisible by $A$. 
We will prove in the lemmas below that for generic $n \in \mathcal{M}$, a sufficient lower bound on the inner sum of \eqref{pw_GCD_bis} holds.

\begin{lem}\label{generic_estimate}
 Let $n = p_1\cdots p_k \in \mathcal{M}$ and $I \in [ \mathcal{P}(n) ]^{j}$ for some $1\leq j \leq k.$  Assume that $T^{2j} < p_I < U^{j/A}$. Then we obtain 
\[\sum_{q \in \mathcal{M}_{n,I}} \frac{1}{\sqrt{q}} \gg_{j,A} \frac{\sqrt{p_I}}{(\log U)^{j/A}}.\]
\end{lem}

\begin{lem}\label{comb_lem}
Let $j > A$. Then we have for $N(1 -o(1))$ many elements $n \in \mathcal{M}$ that
    \begin{equation}\label{number_of_good_PI}|\{I \in [\mathcal{P}(n)]^j, T^{2j} < p_I < U^{j/A}\}| \gg_{j,A} (k^{1/2 - \varepsilon/10})^j.
    \end{equation}
\end{lem}

We now conclude the proof of Theorem \ref{unbounded_GCD}, assuming for now the statements of Lemmas \ref{generic_estimate} and \ref{comb_lem} given, which are proven below.
 By Lemma \ref{comb_lem}, there exists a set $\mathcal{M}' \subseteq \mathcal{M}$ with $|\mathcal{M'}| > N/2$ such that for $n \in \mathcal{M}'$, we have \eqref{number_of_good_PI}.
    Using \eqref{pw_GCD_bis} and Lemma \ref{generic_estimate} we get
    \[\begin{split}\sum_{n \in \mathcal{M}}\sum_{m \in \mathcal{M}} \frac{\gcd(m,n)}{\sqrt{mn}}
   & \gg_{j,A} \sum_{n \in \mathcal{M}'}\sum_{\substack{I \in [\mathcal{P}(n)]^j\\T^{2j}< p_I < U^{j/A} }} \frac{1}{(\log U)^{j/A}} 
   \\ &\gg_{j,A}
    \sum_{n \in \mathcal{M}'}\left(\frac{k^{1/2 - \varepsilon/10}}{(\log U)^{1/A}}\right)^{j}
    \gg N \left(\frac{k^{1/2 - \varepsilon/10}}{(\log U)^{1/A}}\right)^{j}.
    \end{split}
    \]
    Since $\frac{k^{1/2 - \varepsilon/10}}{(\log U)^{1/A}} \asymp (\log N)^{\beta}$
    with $\beta = \beta(\varepsilon) = \varepsilon - \frac{\varepsilon^2}{5} - \varepsilon + \varepsilon^2 > 0$, Theorem \ref{unbounded_GCD} follows by choosing $j > \frac{C}{\beta}$.

\begin{proof}[Proof of Lemma \ref{generic_estimate}]
We aim to find sufficiently many distinct ways of choosing the $j$ prime factors of $q$ such that $q\in \mathcal{M}_{n,I}$. We first choose $q_A:= q_1\cdots q_{j(A-1)/A}$ as a product of $j(A-1)/A$ distinct prime factors of $q$ in the interval $I_A:=(T,T^{A/(A-1)})$, each coprime to $n$.
This can be easily done by induction. Assume that $q_1, \dots, q_{i-1}$ for $1 \leq i \leq j$ have already been chosen to be prime factors of $q_A$. Then, letting $\mathbb{P}_i:= \mathbb{P} \setminus (\mathcal{P}(n)\cup \{q_1,\ldots,q_{i-1}\})$ we can choose $q_i$ as a prime factor of $q_A$ if and only if $q_i \in \mathbb{P}_i.$  The remaining prime factors $q_{(j(A-1)/A) +i}, i=1,\dots, j/A$ will be chosen to not exceed $\left(\frac{p_I}{q_A}\right)^{A/j}$ to ensure that $q \leq p_I$. The squarefreeness and coprimality is ensured again by choosing $q_{(j(A-1)/A) +i} \in \mathbb{P}_{(j(A-1)/A)+i}$.  Noting that every integer $q$ constructed in this way is counted at most $j!$ times, we have the lower bound
     \begin{equation}
         \label{eq:lots_sums}
     \begin{split}\sum_{q \in \mathcal{M}_{n,I}} \frac{1}{\sqrt{q}}
     &\geq \frac{1}{j!}\sum_{\substack{q_1 \in I_A\\ q_1 \in \mathbb{P} \setminus \mathcal{P}(n)}}\frac{1}{\sqrt{q_1}}\sum_{\substack{q_2 \in I_A\\ q_2 \in \mathbb{P}_2}}\frac{1}{\sqrt{q_2}}\ldots \ldots \!\!\!\!\!\!\!\!\!\!\!\!\!\sum_{\substack{q_{j(A-1)/A} \in I_A \\ q_{j(A-1)/A} \in \mathbb{P}_{j(A-1)/A}}}\!\!\!\!\!\!\!\!\!\frac{1}{\sqrt{q_{j(A-1)/A}}}\\&\times\!\!\!\!\!\!\!\!\!\!\
     \sum_{\substack{T < q_{j(A-1)/A+1} < \min\left\{\left(\frac{p_I}{q_A}\right)^{A/j},U\right\}\\ q_{(j(A-1)/A)+1} \in \mathbb{P}_{j(A-1)/A+1}}} \frac{1}{\sqrt{q_{(j(A-1)/A)+1}}}\cdots \!\!\!\!\!\!\!\!\!\!
       \sum_{\substack{T < q_{j} < \min\left\{\left(\frac{p_I}{q_A}\right)^{A/j},U\right\}\\ q_{j} \in \mathbb{P}_j}} \frac{1}{\sqrt{q_{j}}}.
     \end{split}
     \end{equation}
   By the assumption $p_I>T^{2j}$, we have 
    \begin{equation}\label{low_prod}\left(\frac{p_I}{q_1\cdots q_{j(A-1)/A}}\right)^{A/j} = \left(\frac{p_I}{q_A}\right)^{A/j} > \left(\frac{T^{2j}}{T^j}\right)^{A/j} = T^A > T^2.     \end{equation}
        Note that by the assumption $p_I < U^{j/A}$, we have \[\left(\frac{p_I}{q_A}\right)^{A/j} < p_I^{A/j} < U.\]
    This implies that
        \[\begin{split} \sum_{\substack{T < q_{(j(A-1)/A)+1} < \min\left\{\left(\frac{p_I}{q_A}\right)^{A/j},U\right\}}} \ \frac{1}{\sqrt{q_{(j(A-1)/A)+1}}} &\gg 
        \left(\frac{p_I}{q_A}\right)^{A/2j} \frac{1}{\log p_I}
        \gg_{j,A} \left(\frac{p_I}{q_A}\right)^{A/2j} \frac{1}{\log U}.
        \end{split}
        \]
Using that $n$ has at most $k \leq \log N$ many prime divisors and noting that $\sum_{q \leq \log N} \frac{1}{\sqrt{q}} \ll \sqrt{\log N} \leq T^{\varepsilon}$, this implies by \eqref{low_prod} that the removal of these primes from the summation is negligible, which gives
\[\begin{split} \sum_{\substack{T < q_{(j(A-1)/A)+1} < \min\left\{\left(\frac{p_I}{q_A}\right)^{A/j},U\right\}\\         
\\ q_{(j(A-1)/A)+1} \in \mathbb{P}_{j(A-1)/A)+1}}} \ \frac{1}{\sqrt{q_{(j(A-1)/A)+1}}} 
        &\gg_j \left(\frac{p_I}{q_A}\right)^{A/2j} \frac{1}{\log U},
        \end{split}
        \]  and the same bound holds for all variables $q_{(j(A-1)/A)+ i}, i = 1,\ldots, j/A$. Summing over these $j/A$ many inner variables of \eqref{eq:lots_sums}, we obtain the lower bound
        \[\gg_j \sqrt{\frac{p_I}{q_1\cdots q_{j(A-1)/A}}} \frac{1}{(\log U)^{j/A}}.\]
        We are left to compute the $j(A-1)/A$ outer sums of \eqref{eq:lots_sums}, which by Mertens' estimate (again, the prime divisors of $n$ are negligible) implies
        \[\left(\sum_{T < p < T^{1 + 1/(A-1)}}\frac{1}{p}\right)^{j(A-1)/A} \gg_{j,A} 1.\]
        The result follows.
\end{proof}

\begin{proof}[Proof of Lemma \ref{comb_lem}]
We first show that a typical $n \in \mathcal{M}$ should have many large prime divisors. Let $\rho = 1/2 - \varepsilon/10$ and assume that, for $n \in \mathcal{M}$, we have
\begin{equation}\label{nb_large_div}   \lvert\mathcal{P}(n) \cap [T^2,U] \rvert <k^{\rho}.\end{equation} Recalling \eqref{defU}, we get  \[n \leq U^{(k^{\rho})}(T^2)^k < \exp((\log N)^{(1-\varepsilon + 2\varepsilon\rho)})
    \exp(2(\log N)^{1 - 3\varepsilon + 2\varepsilon})) < \exp((\log N)^{1 - \frac{\varepsilon^2}{10}}) = o(N).\]   So it suffices to show that every $n \in \mathcal{M}$ such that $\lvert\mathcal{P}(n) \cap [T^2,U] \rvert  \geq k^{\rho}$ satisfies \eqref{number_of_good_PI}.

    We distinguish two cases depending on the number of prime divisors of $n$ of moderate size.\\

    \underline{Case 1}: $| \mathcal{P}(n) \cap [T^2,U^{1/A}]| \geq \frac{k^{\rho}}{2}$.\\
    Every $p_I$ that arises as a product of a $j$-tuple of primes out of this set satisfies \[T^{2j} < p_I < U^{j/A}.\] Recalling that $j$ does not depend on $N$, for each such $n$, there are at least \[\binom{ \frac{1}{2}k^{\rho}}{j} \gg_j (k^{\rho})^j\] many such $p_I$'s.\\

    \underline{Case 2}: $| \mathcal{P}(n) \cap [T^2,U^{1/A}]| < \frac{k^{\rho}}{2}.$\\ 

    We now construct elements $p_I$ as a product of $j-1$ small prime factors and one very large prime factor. 
We claim that:

\begin{itemize}
    \item[(i)] $\exists p \mid n: p \in [U^{1/A},U]$,
    \item[(ii)] $|\mathcal{P}(n) \cap [T,T^2]|\geq \frac{k}{3}$.
\end{itemize}
    
    Since 
    $|\mathcal{P}(n) \cap  [T^2,U]| \geq k^{\rho}$, (i) follows easily from  
    $$|\mathcal{P}(n) \cap [U^{1/A},U]| \geq k^{\rho} - \frac{k^{\rho}}{2} =\frac{1}{2}k^{1/2 - 2\varepsilon} > 0.$$ 
    
    For (ii), it suffices to prove that
        \begin{equation}\label{eq:not_many_U}|\mathcal{P}(n) \cap [U^{1/A},U]| < k/2.\end{equation}
        Indeed, since we are in Case 2, this implies 
        \[|\mathcal{P}(n) \cap  [T,T^2]| > k-\left(\frac{k}{2} +\frac{k^{\rho}}{2}\right) \geq \frac{k}{3}.\]
Let us assume that \eqref{eq:not_many_U} does not hold. We get \[n > U^{k/2A} > \exp\left(\frac{1}{4A}(\log N)^{1 + \varepsilon}\right)
    > \exp((\log N)^{1 + 2\varepsilon/3}).\] Letting $V := \exp((\log N)^{1 - 3\varepsilon/2})$, we see that any number $m$ that consists of $k$ prime divisors in $[T,V]$ satisfies \[m < V^k < \exp((\log N)^{1 + 2\varepsilon/3}) < n.\]
    Note however that there are $\binom{|\mathbb{P} \cap [T,V]|}{k} > N$ such squarefree numbers. By the minimality property of $\mathcal{M}$, any such $m<n$ would lie in the set $\mathcal{M}$. This contradicts $\vert \mathcal{M}\vert = N$, thus \eqref{eq:not_many_U} follows, which implies (ii).
    
    We now prove that (i) and (ii) allows us to construct sufficiently many good $p_I$'s.
    Indeed, observe that for fixed $\hat{p} \in [U^{1/A},U]$, and
    any $p_{i_1},\ldots,p_{i_{j-1}} \in [T,T^2]$, we have for $p_I := \hat{p}\cdot p_{i_1}\cdots p_{i_{j-1}}$ that
$U^{1/A} < p_I < U T^{2j} $. For $j > A$ and recalling \eqref{defU} this implies
\[T^{2j} < p_I < U^{j/A}.\]
Hence in Case 2 there are at least
     \[\binom{k/3}{j-1} \gg_j k^{j-1}\]
     many such $p_I$'s, concluding the proof.
\end{proof}

\section{Proofs of results on asymptotic approximation}\label{sec:asymptotic}

\subsection{Prerequisites}

We need the following application of the Cauchy--Schwarz inequality.\footnote{It might well be that this inequality is known in the literature, but since we have not found a proper reference, we present its short proof.}
\begin{lem}\label{minimum}
Let $U,$ $V$ be finite sets and $(\alpha_u)_{u\in U}$, $(\beta_v)_{v\in V}$ be non-negative reals. Then
\begin{equation}\label{minimum_basic}
 \sum_{u\in U}\sum_{v\in V}\min\{\alpha_u,\beta_v\}
 \leq
 \left(\sum_{u,u'\in U}\min\{\alpha_u,\alpha_{u'}\}\right)^{1/2}
 \left(\sum_{v,v'\in V}\min\{\beta_v,\beta_{v'}\}\right)^{1/2}.
\end{equation}
If, in addition, $0<\Theta\leq1$ and $\alpha_u\leq\Theta \beta_v$
for every $u\in U$, $v\in V$, then
\begin{equation}\label{minimum_theta}
 \sum_{u\in U}\sum_{v\in V}\min\{\alpha_u,\beta_v\}
 \leq\sqrt{\Theta}
 \left(\sum_{u,u'\in U}\min\{\alpha_u,\alpha_{u'}\}\right)^{1/2}
 \left(\sum_{v,v'\in V}\min\{\beta_v,\beta_{v'}\}\right)^{1/2}.
\end{equation}
\end{lem}

\begin{proof}
We use ``Cauchy--Schwarz on the minimum''. Indeed,
\begin{align*}
 \sum_{u\in U}\sum_{v\in V}\min\{\alpha_u,\beta_v\}
 &=\int_0^\infty
   \left(\sum_{u\in U}\mathds{1}_{[t<\alpha_u]}\right)
   \left(\sum_{v\in V}\mathds{1}_{[t<\beta_v]}\right)\,\mathrm{d}t\\
 &\leq
   \left(\int_0^\infty
     \left(\sum_{u\in U}\mathds{1}_{[t<\alpha_u]}\right)^2\mathrm{d}t
   \right)^{1/2}
   \left(\int_0^\infty
     \left(\sum_{v\in V}\mathds{1}_{[t<\beta_v]}\right)^2\mathrm{d}t
   \right)^{1/2}\\
 &=\left(\sum_{u,u'\in U}\min\{\alpha_u,\alpha_{u'}\}\right)^{1/2}
   \left(\sum_{v,v'\in V}\min\{\beta_v,\beta_{v'}\}\right)^{1/2}.
\end{align*}
Under the additional assumption $\alpha_u\leq\Theta \beta_v$, we have
$\min\{\alpha_u,\beta_v\}=\min\{\alpha_u,\Theta \beta_v\}$. Consequently,
\begin{align*}
 \sum_{u\in U}\sum_{v\in V}\min\{\alpha_u,\beta_v\}
 &=\int_0^\infty
   \left(\sum_{u\in U}\mathds{1}_{[t<\alpha_u]}\right)
   \left(\sum_{v\in V}\mathds{1}_{[t<\Theta \beta_v]}\right)\,\mathrm{d}t\\
 &\leq
   \left(\int_0^\infty
     \left(\sum_{u\in U}\mathds{1}_{[t<\alpha_u]}\right)^2\mathrm{d}t
   \right)^{1/2}
   \left(\int_0^\infty
     \left(\sum_{v\in V}\mathds{1}_{[t<\Theta \beta_v]}\right)^2\mathrm{d}t
   \right)^{1/2}\\
 &=\sqrt{\Theta}
   \left(\sum_{u,u'\in U}\min\{\alpha_u,\alpha_{u'}\}\right)^{1/2}
   \left(\sum_{v,v'\in V}\min\{\beta_v,\beta_{v'}\}\right)^{1/2}.
\end{align*}
\end{proof}
As a second prerequisite, we show that we may assume in Theorem \ref{gen_asymp} that $\psi$ has polynomial decay.

\begin{lem}\label{wlog_psi}
It suffices to prove Theorem \ref{gen_asymp} under the additional assumption that
\begin{equation}\label{wlog_assum}n^{-2}\leq{\psi}(n)\leq n^{-1/2}.
\end{equation}
\end{lem}

\begin{proof} 
Let $\psi:\N\to[0,1/2)$ be monotonically decreasing, $T: \mathbb{N} \to [1, +\infty)$ increasing, and assume \eqref{doubling} and \eqref{extra_div_thm}.
We introduce the functions
\[
\psi_1(n) := \min\{\psi(n),n^{-1/2}\}, \quad 
 \tilde{\psi}(n):=\max\{n^{-2},\psi_1(n)\}, \quad
 \tilde{T}(n) := \min\{T(n),n^2\}.
\]
The statement follows immediately once we have proven the following properties:

\begin{enumerate}[(i)]
    \item $\tilde{\psi}$ is decreasing.
    \item $n^{-2}\leq\tilde{\psi}(n)\leq n^{-1/2}$.
    \item $\tilde{T}$ is increasing.
    \item $\tilde{T}$ satisfies \eqref{doubling}.
    \item $\tilde{T}$ satisfies \eqref{GCD_rough_bound}.
    \item $\sum_{n \in \mathbb{N}}\frac{\tilde{\psi}(n)}{\log \tilde{T}(n)}= \infty$.
    \item $\lambda(W(\tilde{\psi},\A))
 \leq\lambda(W(\psi,\A)).$
\end{enumerate}

Properties (i) - (iv) are obvious from the definition. 
For (v), observe that if $n,m \in \mathcal{M}$ with $n < m$ are both $N^2$-rough, then their prime factorizations differ by at least one prime $p > N^2$. In particular, this implies
$\frac{\gcd(n,m)}{\sqrt{nm}} \leq \frac{1}{N}$, hence the GCD sum is trivially  bounded.

For (vi), we distinguish two cases: if $\psi(n)\leq n^{-1/2}$ for all sufficiently large $n$,
then $\tilde{\psi}(n)\geq\psi(n)$ for all such $n$.
Since $\tilde{T}\leq T$, we obtain
\[
 \sum_{n\in\N}\frac{\tilde{\psi}(n)}{\log\tilde{T}(n)}
 \gg \sum_{n\in\N}
 \frac{\psi(n)}{\log T(n)}
 =\infty.
\]
If $\psi(N)>N^{-1/2}$ for infinitely many $N$, then monotonicity implies for all $N/2\leq n\leq N$ that $\tilde{\psi}(n)\geq N^{-1/2}$.
Using $\tilde{T}(n)\leq n^2$, we therefore obtain
\[
 \sum_{N/2\leq n\leq N}
 \frac{\tilde{\psi}(n)}{\log\tilde{T}(n)}
 \gg \frac{\sqrt N}{\log N}
 \longrightarrow\infty.
\]
Finally, we prove (vii). Note that $\tilde{\psi}(n)\neq\psi_1(n)$ only when
$\psi_1(n)<n^{-2}$, in which case
$\tilde{\psi}(n)=n^{-2}$. Since
\[
 \sum_{n\in\N}
 \lambda\bigl(\{\alpha\in[0,1):
             \|a_n\alpha\|\leq n^{-2}\}\bigr)
 \ll\sum_{n\in\N}n^{-2}<\infty,
\]
the first Borel--Cantelli lemma and ${\psi_1} \leq \psi$ give
\[
 \lambda\bigl(W(\tilde{\psi},\A)\bigr)
 =\lambda\bigl(W(\psi_1,\A)\bigr)
 \leq\lambda\bigl(W(\psi,\A)\bigr).
\]
\end{proof}

\subsection{Proof of Theorem \ref{gen_asymp}}
Let $\psi:\N\to[0,1/2)$ be decreasing, $ T: \mathbb{N} \to (1,+\infty)$ increasing, and assume that \eqref{GCD_rough_bound}, \eqref{doubling} and \eqref{extra_div_thm} hold, as well as, justified by Lemma \ref{wlog_psi}, 
 \eqref{wlog_assum}.
Define $N_i:=2^{(4^i)}$ for $i\geq0$, and put
\[
 \Delta_i:=[N_{i-1},N_i),\qquad
 T_i:=T(N_{i}).
\]
For $n\in\Delta_i$, write
$a_n=\ans\cdot \anl$, with $\ans$ the $T_i$-smooth part, and $\anl$ the
$T_i$-rough part. Define
\[
A_n := \bigcup_{\substack{0\leq a<a_n}}
 \left(\frac a{a_n}-\frac{\psi(n)}{a_n},
       \frac a{a_n}+\frac{\psi(n)}{a_n}\right), \quad
 A_n^*:=\bigcup_{\substack{0\leq a<a_n\\\gcd(a,\ans)=1}}
 \left(\frac a{a_n}-\frac{\psi(n)}{a_n},
       \frac a{a_n}+\frac{\psi(n)}{a_n}\right).
\] Note that we have 
\begin{equation}\label{meas_starsets}
    \lambda(A_n^*)= 2\psi(n) \frac{\varphi(\ans)}{\ans} \gg \psi(n) \prod_{p \leq T_i}\left(1 - \frac{1}{p}\right) \gg \frac{\psi(n)}{\log T_i} \gg \frac{\psi(n)}{\log T(n)},
\end{equation} 
using the doubling assumption \eqref{doubling} in the last inequality. Thus, using \eqref{extra_div_thm}, we obtain
\begin{equation}\label{A_n*_div}\sum_{n = 1}^{+\infty} \lambda(A_n^*) = \infty.\end{equation}
We now partition $A_n^*$ in a set of \textit{reduced} fractions, i.e.,
writing

\begin{equation}\label{eq:Ad'}
A_d'(\rho) := \bigcup_{\substack{0 \leq a < d\\\gcd(a,d) = 1}} \left[\frac{a}{d} - \rho, \frac{a}{d} + \rho\right],\end{equation}
we have

\begin{equation}\label{Aq_part_} A^{*}_n := \mathop{\dot{\bigcup}}_{d\mid^* a_n} A_{d}'(\rho_n),\end{equation}
where
\begin{equation}\label{def_mid*}\rho_n := \frac{\psi(n)}{a_n}, \quad d \mid^* a_n :\Leftrightarrow d = e\ans\quad \text{ for some }\quad e \mid \anl.\end{equation} 
Note that the condition $\mid^*$ depends on the block $\Delta_i$ where $n$ belongs to. 
We further let 
\begin{equation}\label{defSd} S_d := \{n: d \ms a_n\},\end{equation}
and
\[\psi_{\A}^*(d) := \begin{cases}
    d\sup_{n: d \ms a_n}\rho_n &\text{ if } S_d \neq \emptyset, \\
    0 &\text{ if } S_d = \emptyset.
\end{cases}\] 
We claim that if
\begin{equation}\label{catlin_div}\sum_{d \in \N}\frac{\varphi(d)}{d}\psi_{\A}^*(d) = \infty,\end{equation}
then $\lambda(\limsup_{n \to \infty}A_n^*) = 1$ (and thus in particular, $\lambda(\limsup_{n \to \infty}A_n) = 1$) follows. Indeed, by the Koukoulopoulos--Maynard theorem \cite[Theorem 1]{KM20}, for almost every $\alpha$, there exist infinitely many solutions $(a,d) \in \mathbb{Z} \times \mathbb{N}$ to
\begin{equation}\label{reduced_sol_proof}\left\lvert \alpha - \frac{a}{d}\right \rvert \leq \frac{\psi_{\A}^*(d)}{d}, \quad \gcd(a,d) = 1.\end{equation}
 Since $\rho_n = \psi(n)/a_n\to0$, the supremum is attained at some $n\in S_d$.
For this $n$ we have $\psi(n)/a_n=\psi_{\A}^*(d)/d$ and hence
$A_d'(\psi_{\A}^*(d)/d)\subseteq A_n^*$. Each fixed $n$ has only finitely many
divisors $d\ms a_n$, so infinitely many $d$ with \eqref{reduced_sol_proof} imply infinitely many $n$ with $\alpha \in A_n^*$.
Hence, $\lambda(\limsup_{n \to \infty}A_n) = 1$ follows from \eqref{catlin_div}, which we are left to establish.

For $J \in \mathbb{N}$, we write $D_J=a_{N_J-1}$, and aim to show

\begin{equation}\label{catlin_estimate_}
 \sum_{d\leq D_J}\frac{\varphi(d)}d\psi_{\A}^*(d)
 \gg \sum_{2\leq n<N_J}\lambda(A_n^*),
\end{equation}
with an implied constant that is uniform in $J$.  
Letting $J\to\infty$ and using \eqref{A_n*_div} then shows \eqref{catlin_div}.
 Since $d\ms a_n$ implies $d\leq a_n$, we have by \eqref{Aq_part_}
\begin{equation}\label{inclusion_reduced}
 \bigcup_{1\leq n<N_J}A_n^*
 \subseteq\bigcup_{d\leq D_J}A_d'(\psi_{\A}^*(d)/d).
\end{equation}
Hence, we obtain
\begin{equation}\label{cs1_}
 \sum_{1\leq n<N_J}\lambda(A_n^*) =  2\sum_{1\leq n<N_J}\sum_{ d \ms a_n} \varphi(d) \rho_n
 \ll \sum_{d\leq D_J}\varphi(d)
 \sum_{\substack{n\in S_d\\n<N_J}}\rho_n.
\end{equation} 
Now fix $d \leq D_J, 1 \leq i,j \leq J$
and define
\begin{align*}
 C_{i,j}(d)&
 :=\varphi(d)\sum_{n\in S_i(d)}\sum_{m\in S_j(d)}
 \min\{\rho_n,\rho_m\},\\
 C_J(d)&:=\sum_{1\leq j\leq i\leq J}C_{i,j}(d) \asymp \varphi(d) \sum_{\substack{n,m\in S_d\\n,m<N_J}}\min\{\rho_n,\rho_m\},
\end{align*} where 
\begin{equation}\label{def_Si} S_i(d)=S_d\cap\Delta_i.\end{equation}
Since by definition, $\rho_n \leq \psi_{\A}^*(d)/d$ for all $n \in S_d$, we get
\[
 \left(\sum_{\substack{n\in S_d\\n<N_J}}\rho_n\right)^2
 \leq\frac{\psi_{\A}^*(d)}d\sum_{\substack{n,m\in S_d\\n,m<N_J}}
 \min\{\rho_n,\rho_m\}
 \ll\frac{\psi_{\A}^*(d)C_J(d)}{d\varphi(d)}.
\]
Hence, by \eqref{cs1_} and Cauchy--Schwarz we get 
\begin{equation}\label{cs2}
 \left(\sum_{1\leq n<N_J}\lambda(A_n^*)\right)^2
 \ll \left(\sum_{d\leq D_J}
 \sqrt{\frac{\varphi(d)}d\psi_{\A}^*(d)}\sqrt{C_J(d)}\right)^2\\
 \leq\left(\sum_{d\leq D_J}\frac{\varphi(d)}d\psi_{\A}^*(d)\right)
             \left(\sum_{d \leq D_J} C_J(d)\right).
\end{equation} 
Therefore, in order to prove \eqref{catlin_estimate_}, it suffices to prove
\begin{equation}\label{corr_bound}
 \sum_{d \leq D_J} C_J(d) \leq 2  \sum_{d \leq D_J}\sum_{1 \leq i \leq j \leq J} C_{i,j}(d) \ll\sum_{2\leq n<N_J}\lambda(A_n^*).
\end{equation} 
Recall that by \eqref{wlog_assum}, we have $\frac{1}{n^{2}}\leq\psi(n)\leq \frac{1}{\sqrt{n}}$. This implies for $m \in \Delta_j, n \in \Delta_i$ that
\[
 \frac{\rho_n}{\rho_m}
 =\frac{\psi(n)}{\psi(m)}\frac{a_m}{a_n}
 \leq\frac{m^2}{n^{1/2}}\leq\frac{N_j^2}{\sqrt{N_{i-1}}},
\]
where we used that $(a_n)_n$ is increasing.
We apply Lemma \ref{minimum}: for $j \in \{i-1,i\}$ in the form of \eqref{minimum_basic}, for $j \leq i-2$ in the form of \eqref{minimum_theta}, with $\Theta=N_j^2/\sqrt{N_{i-1}}$.
By the arithmetic–geometric mean inequality, as well as recalling $N_j = 2^{(4^j)}$, we obtain
\begin{equation}\label{cross}
 C_{i,j}(d) \ll  2^{-(i-j)}\left(C_{i,i} + C_{j,j}\right).
\end{equation}
Consequently, we get
\[C_J(d) \ll \sum_{1\leq j\leq i \leq J} C_{i,j} \ll \sum_{i=1}^J C_{i,i}(d),\]
thus it suffices to control the terms $C_{i,i}$.
For fixed $i$, write $d=\ds\cdot \dl$ for its $T_i$-smooth and rough parts. 
Note that $n,m\in S_i(d)$ implies
$\ds=\ans=\ams$ and $\dl\mid(\anl,\aml)$. Thus
\begin{align*}
 \sum_{d \leq D_J} C_{i,i}(d)
 &\ll\sum_{\ds}\varphi(\ds)
 \sum_{\substack{n,m\in\Delta_i\\\ans=\ams=\ds}}
 \sum_{\dl\mid(\anl,\aml)}\varphi(\dl)\min\{\rho_n,\rho_m\}\\
 &\ll\sum_{\ds}\frac{\varphi(\ds)}{\ds}
 \sum_{\substack{n,m\in\Delta_i\\ m \leq n \\ \ans=\ams=\ds}}
 (\anl,\aml)\min\left\{\frac{\psi(n)}{\anl},
                         \frac{\psi(m)}{\aml}\right\}.
\end{align*}
For $m\leq n$, we have $\psi(n)\leq\psi(m)$. Moreover, $ \ans = \ams$, so since $(a_n)_n$ is increasing, $\anl\geq\aml$. It follows that
\[
 \min\left\{\frac{\psi(n)}{\anl},\frac{\psi(m)}{\aml}\right\}
 \leq\frac{\min\{\psi(n),\psi(m)\}}{\sqrt{\anl\aml}}.
\]
Note that the integers $\anl$ are all distinct once $\ans= d$ is fixed, and the set
\[\mathcal{M}_{i,d} = \{\anl: n \in S_i(d)\}\] has cardinality at most $N_i$.
Thus applying Theorem \ref{doublelog_thm} to the inner sum shows
\begin{align}\label{diag}
 \sum_d C_{i,i}(d)
& \ll\sum_{\ds}\frac{\varphi(\ds)}{\ds}\sum_{\substack{n,m\in\Delta_i\\m \leq n\\\ans=\ams=\ds}}
 \frac{(\anl,\aml)}{\sqrt{\anl\aml}}\min\left\{\psi(n),\psi(m)\right\} \nonumber\\ 
  & \ll  \sum_{\ds}\frac{\varphi(\ds)}{\ds}   \sum_{\substack{n\in\Delta_i\\\ans=\ds}}\psi(n) = \sum_{n\in\Delta_i} \psi(n)\frac{\varphi(\ans)}{\ans}  
 \ll \sum_{n\in\Delta_i}\lambda(A_n^*).
\end{align}
Summing over $i$ shows \eqref{corr_bound}, which proves Theorem \ref{gen_asymp}.

\subsection{First part of Theorem \ref{nonincr_thm}}
As above, let $N_i := 2^{(4^i)}$, now choosing $T_i := N_i^2$ and $\psi(n):=\tfrac{1}{n}.$
Again, it suffices to prove \eqref{catlin_div}. To this end, define
\[\psi_i^*(d) := \begin{cases}
    d\sup_{n \in \Delta_i: d \ms a_n}\rho_n &\text{ if } S_i(d) \neq \emptyset,\\
        0 &\text{ if } S_i(d) = \emptyset,
\end{cases}\]  where $S_i(d)$ is defined in \eqref{def_Si}.
We note that it follows immediately that $\psi_{\A}^*(d) = \sup_{i}\psi_i^*(d)$
and for fixed $d$, $\lim_{i \to \infty} \psi_i^*(d) = 0$.
Assuming by contradiction that
\begin{equation}\label{eq:contra_conv}\sum_{d}\frac{\varphi(d)}{d}\psi_{\A}^*(d) < \infty\end{equation}
implies that for every $\varepsilon > 0$, there exists $X = X(\varepsilon)$ such that
\[\sum_{d \geq X}\frac{\varphi(d)}{d}\psi_{\A}^*(d) < \varepsilon.\]
However, since $\psi_{\A}^*(d) = \sup_{i}\psi_i^*(d)$, this implies that for every $i \in\mathbb{N}$,
\[\sum_{d \geq X}\frac{\varphi(d)}{d}\psi_i^*(d) < \varepsilon.\]
Since $\lim_{i \to \infty} \psi_i^*(d) = 0$, for $i$ sufficiently large, we also have
\[\sum_{d <X}\frac{\varphi(d)}{d}\psi_i^*(d) < \varepsilon,\]
and consequently, 
\[\sum_{d \in \mathbb{N}}\frac{\varphi(d)}{d}\psi_i^*(d) < 2\varepsilon.\]
Hence we can contradict the assumption of \eqref{eq:contra_conv} if we can show that for all $i \in \mathbb{N},$
\begin{equation}\label{mass_gg1}M_i := \sum_{d}\frac{\varphi(d)}{d}\psi_{i}^*(d) \gg 1.\end{equation}
Note that in place of \eqref{meas_starsets}, we get
$\lambda(A_n^*) \gg \frac{\psi(n)}{\log n}$, and hence,
\[\sum_{n \in \Delta_i}\lambda(A_n^*) \gg 1.\]
Consequently, it suffices to show (compare with \eqref{catlin_estimate_})
\begin{equation}\label{catlin_estimate_Mi}
 \sum_{d}\frac{\varphi(d)}{d}\psi_{i}^*(d)
 \gg \sum_{n \in \Delta_i}\lambda(A_n^*).
\end{equation}
Mirroring the proof of Theorem \ref{gen_asymp}, \eqref{catlin_estimate_Mi} reduces to showing

 \[\sum_d C_{i,i}(d) 
 =\sum_d\varphi(d)\sum_{n,m\in S_i(d)}
 \min\{\rho_n,\rho_m\} 
 \ll \sum_{n \in \Delta_i}\lambda(A_n^*).
 \]
 We note that here, in contrast to the proof of Theorem \ref{gen_asymp}, no ``non-diagonal terms'' $C_{i,j}, i \neq j$ arise.
 Following the proof further, we have
 \[\sum_d C_{i,i}(d)
 \ll 
 \sum_{\ds}\frac{\varphi(\ds)}{\ds}
 \sum_{\substack{n,m\in\Delta_i\\ \anl \leq \aml \\ \ans=\ams=\ds}}
 (\anl,\aml)\min\left\{\frac{\psi(n)}{\anl},
                         \frac{\psi(m)}{\aml}\right\}.
  \]
 Clearly, the terms where $\anl = \aml, \ans=\ams$ only contribute
 \[\sum_{d }\frac{\varphi(\ds)}{\ds}\sum_{\substack{m\in\Delta_i\\ \ams=\ds}} \psi(m) \ll \sum_{n \in \Delta_i}\lambda(A_n^*).\] 
Note that by $T_i = N_i^2$-roughness,  for $\anl \neq \aml$,
\[ \gcd(\anl,\aml) \min\left\{\frac{\psi(n)}{\anl},
                         \frac{\psi(m)}{\aml}\right\} \leq  \gcd(\anl,\aml)  \frac{\sqrt{\psi(m)\psi(n)}}{\sqrt{\anl \aml}} \leq \frac{1}{N_i}  \sqrt{\psi(m)\psi(n)}.\]
Hence, by Cauchy--Schwarz and $|\Delta_i| \leq N_i$,

\[\begin{split} \sum_{\ds}\frac{\varphi(\ds)}{\ds}
 \sum_{\substack{n,m\in\Delta_i\\ \ans=\ams=\ds}}
 (\anl,\aml)\min\left\{\frac{\psi(n)}{\anl},
                         \frac{\psi(m)}{\aml}\right\}
                         &\leq \frac{1}{N_i} \sum_{\ds}\frac{\varphi(\ds)}{\ds}\left(\sum_{\substack{m\in\Delta_i \\\ams=\ds}}\sqrt{\psi(m)}\right)^2  
                         \\& \ll  \sum_{\ds}\frac{\varphi(\ds)}{\ds}\sum_{\substack{m\in\Delta_i \\\ams=\ds}} \psi(m) \\&\ll  \sum_{n \in \Delta_i}\lambda(A_n^*).
                         \end{split}\] 
This proves \eqref{catlin_estimate_Mi}, and shows that
for almost every $\alpha$,
\[\lVert a_n\alpha\rVert \leq \frac{1}{n} \text{ for i.m. }n \in \mathbb{N}.\]
Note that by Cassels's constant invariance principle \cite[Lemma 9]{Cassels}, we have $\lambda(W({\A},\psi)) = \lambda(W({\A},\varepsilon\psi))$
    for every $\varepsilon > 0$. Hence the first part of Theorem \ref{nonincr_thm} follows by $\varepsilon \to 0$, and the continuity of measures.
    \subsection{Second part of Theorem \ref{nonincr_thm}: The counterexample}\label{sec:counterex}

    In this subsection we show the second part of Theorem \ref{nonincr_thm}, i.e.:
 for every $\psi$ with $\psi(n) = o(1/n)$, there exists a sequence of (in general not increasing) distinct integers $\A = (a_n)_n$ such that
    \[\lambda(W({\A},\psi)) = 0.\]
As a first reduction, note that we may assume that $n \mapsto n \psi(n)$ is monotonically decreasing, by replacing $\psi$ with the majorant $\hat{\psi}(n) := \frac{\sup_{m \geq n} m \psi(m)}{n}$, which still satisfies $\hat{\psi}(n) = o(1/n)$.
Write $\Delta_i := [2^i,2^{i+1})$. We may further
restrict attention to those $\psi$ such that $\psi(n) \equiv \psi_k, n \in \Delta_k$ with $\frac{\psi_k}{\psi_{k+1}} \geq 2$. This follows from choosing the further majorant $\tilde{\psi}(n) = \hat{\psi}(2^k) =: \psi_k, n \in \Delta_k$, and changing at most finitely many values for small $k$.

    \begin{lem}\label{counterex_block}
    Let $(\psi_k)_k$ be as above, $\varepsilon > 0$ be arbitrary, and $k = k(\varepsilon,\psi)$ be such that $\psi_k \leq \frac{\varepsilon}{2^k}$. Then we can find $R = R(k,\varepsilon)$ such that $\psi_{k+R} \leq \frac{\varepsilon/2}{2^{k+R}}$ and a set of distinct integers
    $(a_n)_{n = 2^k}^{2^{k+R}-1}$ with $\min_{n}a_n$ arbitrarily large, such that
    \[\lambda\left(\bigcup_{k \leq i < k+R} \bigcup_{n \in \Delta_i} A_n\right)
    \ll \varepsilon.
    \]
    \end{lem}

 Once Lemma \ref{counterex_block} is proven, the second part of Theorem \ref{nonincr_thm} follows immediately: 
    choose $\varepsilon_j = 2^{-j}, j \in \mathbb{N}$. We may then apply Lemma \ref{counterex_block} for $j =1$: we find $k_1,R_1$ such that for  
    $M_1 := \{2^{k_1} \leq n < 2^{k_1 + R_1}\}$, we have
    \[\lambda\left(\bigcup_{n \in M_1} A_n\right)
    \ll \varepsilon_1.\]
    Next, note that by construction $\psi_{k_1+R_1} \leq \frac{\varepsilon_1/2}{2^{k_1 + R_1}} = \frac{\varepsilon_2}{2^{k_1 + R_1}} $, thus we can apply Lemma \ref{counterex_block} with $\varepsilon_2$ to
    $k_2 = k_1 + R_1$, and define $R_2,M_2$ analogously such that 
    \[\lambda\left(\bigcup_{n \in M_2} A_n\right)
    \ll \varepsilon_2.\]
    Note that in Lemma \ref{counterex_block}, we can assume $\min_{n \in M_2}a_n$ arbitrarily large, thus in particular, larger than $\max_{n \in M_1}a_n$, which ensures that the sequence $(a_n)_{n \leq M_2}$ contains only distinct integers. Iterating this procedure inductively, we obtain
    \[\lambda(W({\A},\psi))= \lambda(\limsup_{n \to \infty}A_n) \leq
    \lim_{j \to \infty} \sum_{i \geq j} \lambda\left(\bigcup_{n \in M_i} A_n\right)
    \ll \lim_{j \to \infty} \sum_{i \geq j} 2^{-i} = 0.
    \]
    Hence, we are left to prove Lemma \ref{counterex_block}.

\begin{proof}[Proof of Lemma \ref{counterex_block}]
The construction builds on the idea of the classical counterexample by Duffin and Schaeffer \cite{DS41}. We begin by choosing the smooth integer $a = (2^{R+1})!$.
    We then define $S_0(a) := \{a\}$, and for $1 \leq i < R$, the sets
    $S_i(a) := \left\{\frac{a}{d}, d \in \Delta_i\right\}$. These are sets of $2^i$ integers, ensured by choosing $a$ to be the large factorial. We then define
    for $b \in S_i(a)$ that $\chi(b) := \psi_{k+i}$, and $\chi(a) := \psi_k$.
We claim that for
\[B_n(\chi) = \{\alpha \in [0,1): \lVert n \alpha \rVert \leq \chi(n)\}\]
(note that in our notation, $A_n = B_{a_n}$), we have
\[\lambda\left(\bigcup_{0 \leq i < R}\bigcup_{b \in S_i(a)} B_b(\chi)\right)
\ll \lambda(B_a(\chi)),
\]
by showing 
\begin{equation}\label{pw_incl}\bigcup_{0 \leq i < R}\bigcup_{b \in S_i(a)} B_b
\subseteq B_a^{(2)} := \{\alpha \in [0,1): \lVert a \alpha \rVert \leq 2 \chi(a)\}.
\end{equation}
To see this, note that for any $i$  and any $b \in S_i(a)$ we have
$b \mid a$, and in fact, $\frac{a}{b} \asymp 2^i$.
The inclusion $B_b \subseteq B_a^{(2)}$ thus holds if
\[ \frac{\chi(b)}{b} \leq \frac{2\chi(a)}{a}.\]
This is true by construction, since
\[\frac{2b}{a} \geq 2^{-i}, \quad \frac{\chi(b)}{\chi(a)} = \frac{\psi_{k+i}}{\psi_{k}} \leq 2^{-i}.\]
We now construct ``copies'' of the above systems, by exploiting the fact that multiplying the system by arbitrary integers does not change the inclusions whatsoever. Let
\[a_n = p(n) \cdot (2^{R+1})!, \quad n \in \Delta_k,\] where $\{p(n): n \in \Delta_k\}$ is a set of distinct primes all exceeding $(2^{R+1})!$. Consequently, defining $S_{i,n} := p(n) \cdot S_i(a)$, we see that $S_{i,n} \cap S_{j,m} = \emptyset$ for $(i,n) \neq (j,m)$.
We now define
\[\{a_n: n \in \Delta_{k+i}\} = \dot{\bigcup_{n \in \Delta_k}}S_{i,n},\]
which is a set of cardinality $2^{i+k}$. We may pick an arbitrary order of the elements in this set, and assign to each the value $\chi(p(n) b) := \psi_{k+i}$
for all $b \in S_{i}(a), n \in \Delta_k$.
    By \eqref{pw_incl}, this now implies
    \[\bigcup_{0 \leq i < R}\bigcup_{n \in \Delta_{k+i}}A_n
    = \bigcup_{n \in \Delta_k} \bigcup_{0 \leq i < R} \bigcup_{m \in S_{i,n}}B_m(\chi)
    \subseteq \bigcup_{n \in \Delta_k} A_n.
    \]
    Now note that 
    \[\lambda\left(\bigcup_{n \in \Delta_k} A_n\right) \ll
    \sum_{n \in \Delta_k}\lambda(A_n) \ll \sum_{n \in \Delta_k}\psi(n) \ll \sum_{n \in \Delta_k}\frac{\varepsilon}{n} \leq 
    \varepsilon.\]
    This concludes the proof by remarking that choosing the primes $p(n)$ all sufficiently large ensures that $\min_{n} a_n$ is arbitrarily large.
\end{proof}

\subsection{Proof of Corollary \ref{asymp_thm_cases}}
    Let
    \[I := \left\{n \in \N: 
    \frac{\log\log a_n}{\log n\,\log\log n} <\frac{1}{(\log\log n)^{2}}
    \right\}.\]
    We distinguish two cases.
    \underline{Case $1$}: $\sum_{n \in I} \frac{1}{n \log n \log \log n} = \infty$. 
    Then $\sum_{n \in I} \frac{\log \log a_n}{n \log n \log \log n} \frac{\varphi(a_n)}{a_n} = \infty$. So by the Koukoulopoulos-Maynard theorem \cite[Theorem $1$]{KM20}, we get infinitely many solutions with $n \in I$, so in particular infinitely many $n \in \N$.\\
    
  \noindent  \underline{Case $2$}: Let $J = \mathbb{N}\setminus I$, and assume
    $\sum_{n \in J} \frac{1}{n \log n \log \log n} = \infty$. Then let $J = \{j_1 < j_2 \ldots\}$ and consider the subsequence $b_i = a_{{j_i}}$, and set $T(i)=\exp\left(20\frac{\log i}{\log \log i}\right)$ and
    $\psi(i) = \frac{1}{j_i(\log \log j_i)^{2}}$. Then we have
    (since $x\mapsto (\log\log x)/\log x$ is decreasing, and $i \leq j_i$)
    \begin{align*}\sum_{i} \frac{\psi(i)}{\log T(i)}
    \asymp     \sum_{i}\frac{\psi(i) \log \log i}{\log i} 
    \gg \sum_{j \in J}\frac{\log \log j}{j \log j(\log \log j)^{2}}
    \geq \sum_{j \in J}\frac{1}{j \log j \log \log j} = 
    \infty.\end{align*} By Theorem \ref{doublelog_thm}, 
     $T(i)=\exp\left(20\frac{\log i}{(\log \log i)}\right)$ is an admissible choice, hence we can apply Theorem \ref{gen_asymp}, and the result follows.

\section{Proofs of results on uniform approximation}\label{sec:uniform}

 \subsection{Proof of Theorem \ref{main_thm_unif}}

We proceed along similar lines as in the proof of Theorem \ref{nonincr_thm}.
Let $N_k = 2^k$, $\Delta_k = (N_{k-1},N_k]$, let $\psi_k(n) = \frac{k^{2 + \varepsilon}}{10 \cdot 2^k}$ for $n \in \Delta_k$, and set $T_k = N_k^2$.
Suppose that for almost all $\alpha$

\begin{equation}\label{goal_uniform}\mathcal{N}_k(\alpha) := \sum_{N_{k-1} < n \leq N_k}\mathds{1}_{[\lVert a_n \alpha\rVert \leq \psi_k(n)]} \geq 1\end{equation}
for all but finitely many $k$, then this implies the uniform statement
\[\lVert a_n\alpha\rVert \leq \frac{(\log N)^{2 + \varepsilon}}{N}, \quad  1 \leq n\leq N\]
for all but finitely many $N$, as desired. Writing
$a_n = \ans \cdot \anl$ for the $T_k$-smooth and rough parts of $a_n$, where $n \in \Delta_k$, we define

\[
 A_n^*(\psi_k)=\bigcup_{\substack{0\leq a<a_n\\(a,\ans)=1}}
 \left(\frac a{a_n}-\frac{\psi_k(n)}{a_n},
       \frac a{a_n}+\frac{\psi_k(n)}{a_n}\right).
\]  We clearly have 
\[\mathcal{N}_k(\alpha) \geq \mathcal{N}_k^*(\alpha) := \sum_{n \in \Delta_k}\mathds{1}_{[\alpha \in A_n^*(\psi_k)]}. \]
Similar to before, we partition
\[A^{*}_n(\psi_k)
=  \bigcup_{d\mid^* a_n}A_{d}'(\psi_k(n)/a_n),
 \] 
 with $A_d'$ and $\mid^*$ defined as in \eqref{eq:Ad'} and \eqref{def_mid*}.
 Let \[S_k(d):=\{n \in \Delta_k: d \ms a_n\} \qquad \text{ and } \qquad D_k= \{ d \in \mathbb{N}, S_k(d) \neq \emptyset\}.\] 
With these notations we have
 \[\mathcal{N}_k^*(\alpha) = \sum_{n \in \Delta_k} \sum_{d \ms a_n} \mathds{1}_{[\alpha \in A_d'(\psi_k(n)/a_n)]} \geq \sum_{d \in D_k}\mathds{1}_{[\alpha \in A_d'(\psi^*_k(d)/d)]} =: \mathcal{D}_k(\alpha),\]
 
where $\psi^*_k(d) := \begin{cases}
    d\sup\limits_{\substack{2^{k-1} < n \leq 2^k\\ d \ms a_n}}\psi_k(n)/a_n &\text{ if } S_k(d) \neq \emptyset, \\ 
    0 &\text{ if } S_k(d) = \emptyset. 
    \end{cases}$\\
We claim that
\begin{equation}\label{uniform_left_to_show}\Psi_k := \int_{0}^1 \mathcal{D}_k(\alpha) \,\mathrm{d}\alpha \asymp \sum_{d \in D_k}\frac{\varphi(d)\psi^*_k(d)}{d} \gg k^{1+\varepsilon/2}.\end{equation}
As in \eqref{meas_starsets}, we have

\[\lambda(A_n^*(\psi_k)) = 2\psi_k(n) \frac{\varphi(\ans)}{\ans} \gg \frac{\psi_k(n)}{\log T_k}  \gg \frac{\psi_k(n)}{\log n} \gg \frac{k^{1+\varepsilon}}{2^k}.\]
Thus $\sum_{n \in \Delta_k}\lambda(A_n^*(\psi_k)) \gg k^{1+\varepsilon}$,
and it remains to show, as in the non-uniform variant,

\begin{equation}\label{low_Dk} \Psi_k \asymp \sum_{d \in \mathcal{D}_k}\frac{\varphi(d)}d\psi^*_k(d)
 \gg \sum_{n \in \Delta_k}\lambda(A_n^*(\psi_k)).\end{equation} 
Since the function  $T(n) := T_k \geq n^2,n \in \Delta_k$ satisfies 
\eqref{GCD_rough_bound}, following exactly along the final lines of the proof of Theorem \ref{nonincr_thm}, we obtain \eqref{low_Dk}, and hence \eqref{uniform_left_to_show} follows.

We now apply the variance estimate of Koukoulopoulos--Maynard--Yang \cite[Theorem $2$]{KMY24}: for every $\delta > 0$,
\[
\int_{0}^1 \left(\mathcal{D}_k(\alpha) - \Psi_k\right)^2 \,\mathrm{d}\alpha
\ll_{\delta} \Psi_k^{1+\delta}.
\]
Consequently, by Chebyshev's inequality, using \eqref{uniform_left_to_show},
\[\begin{split}\lambda\left(\{\alpha: \mathcal{N}_k(\alpha) = 0\}\right)
&\leq \lambda\left(\{\alpha: \mathcal{D}_k(\alpha) = 0\}\right)
\leq \lambda\left(\left\{\alpha: \left|\mathcal{D}_k(\alpha) - \Psi_k\right| \geq \frac{\Psi_k}{2} \right\}\right) \\&\ll_{\delta} \frac{{\Psi_k}^{1+\delta}}{\Psi_k^ 2}\ll \frac{1}{k^{(1 + \varepsilon/2)(1 - \delta)}}.
\end{split}
\]
Choosing $\delta$ sufficiently small, the above is bounded by
$\frac{1}{k^{(1 + \varepsilon/4)}}$.
Summing over $k$, the convergence Borel--Cantelli lemma implies \eqref{goal_uniform}, which concludes the proof.

\subsection{Proof of Corollary \ref{unif_thm_cases}}
 By Theorem \ref{main_thm_unif}, for Lebesgue almost every $\alpha\in\R$ and all
sufficiently large $N\in\mathbb N$, there exists $1\leq n\leq N$
such that
\[
 \lVert a_n\alpha\rVert
 \leq \frac{(\log N)^{2+\varepsilon}}{N}.
\] Hence it suffices to show that, for Lebesgue almost every $\alpha\in\R$ and all but finitely many $N$, there exists $1\leq n\leq N$
such that
 \[\lVert a_n\alpha\rVert
 \leq \frac{(\log N)^{1+\varepsilon}}{N}\log \log a_N.\] Corollary \ref{unif_thm_cases} then follows by intersecting the two sets of full measure and taking the maximum of the finitely many exceptions. 
 Let $\mathcal{A}_N = \{a_n: 1 \leq n \leq N\}.$
For every $N \in \mathbb{N}$, we define  \[\psi_N(n) := \mathds{1}_{[n \in \mathcal{A}_N]}\frac{(\log N)^{1 + \varepsilon}}{2N} \log \log a_N \qquad \text{ and }\qquad \Psi_N' := \sum_{n \in \mathbb{N}}\frac{2\varphi(n)\psi_N(n)}{n}.\] 
Defining 
\begin{equation}\label{defN}\mathcal{N}(N,\alpha) = \#\{n \in \mathbb{N}: \lVert n\alpha \rVert \leq 2\psi_N(n)\},\end{equation} we will show that for almost every $\alpha$, there exists a $N_0 = N_0(\alpha)$ such that for 
all $N \geq N_0$, we have $\mathcal{N}(N,\alpha) \geq 1$.
Since $\frac{n}{\varphi(n)} \ll \log \log n \leq \log \log a_N$, we have $\frac{\varphi(n)\psi_N(n)}{n} \gg \frac{(\log N)^{1+\varepsilon}}{N}$ for $n \in \mathcal{A}_N$. Thus, 
for $N$ sufficiently large, 
\begin{equation}\label{eq_size_of_Psi}
\Psi_N' \geq (\log N)^{1 + \varepsilon/2}.
\end{equation} Further, let us write \begin{equation}\label{defNprime}\mathcal{N}'(N,\alpha) = \#\{n \in \mathbb{N}: \lVert n\alpha \rVert' \leq \psi_N(n)\} \qquad \text{ where } \qquad \lVert nx\rVert' := \min_{\gcd(a,n) =1}\lvert nx -a \rvert.\end{equation} 
By \cite[Theorem 2]{KMY24}, we have for all $N \in \mathbb{N}$,

\[
\int_{0}^{1} \left(\mathcal{N}'(N,\alpha) - \Psi_N'\right)^2 \,\mathrm{d}\alpha \ll_{\varepsilon}
(\Psi_N')^{1 + \varepsilon/10}.
\] Thus by Chebyshev's inequality we have for all $N \in \mathbb{N}$,

\begin{equation}\label{chebyshev}
\lambda\left(\left\{\alpha \in [0,1): \mathcal{N}'(N,\alpha) \leq \frac{\Psi_N'}{2} \right\}\right) \ll 
\frac{(\Psi_N')^{1+ \varepsilon/10}}{(\Psi_N')^2} = \frac{1}{(\Psi_N')^{1 - \varepsilon/10}}.
\end{equation}

Now, for any $k \in \N$, let $N_k :=  2^k.$
Combining \eqref{eq_size_of_Psi} and \eqref{chebyshev} implies

\[
\sum_{k \in \mathbb{N}}
\lambda\left(\left\{\alpha \in [0,1): \mathcal{N}'(N_k,\alpha) \leq \frac{\Psi_{N_k}'}{2} \right\}\right) 
\ll \sum_{k \in \N} \frac{1}{((\log N_k)^{1 + \varepsilon/2})^{1 - \varepsilon/8}}
\ll \sum_{k \in \N} \frac{1}{k^{(1 - \varepsilon/8)(1 + \varepsilon/2)}}< \infty,
\]
where we assumed $\varepsilon$ to be sufficiently small. Therefore by the first Borel\,--\,Cantelli lemma, for almost every $\alpha$, there exists $k_0(\alpha) \in \mathbb{N}$ such that for all $k \geq k_0(\alpha)$, we have $\mathcal{N}'(N_k,\alpha) \geq \frac{\Psi_{N_k}'}{2} > 0$.
Let $N$ be such that $N_k \leq N < N_{k+1}$, so that $\frac{N}{N_k} \leq 2.$ Thus
\[\begin{split}
\mathcal{N}(N,\alpha) &= \#\{n \in \mathbb{N}: \lVert n\alpha \rVert \leq 2\psi_N(n)\}
\\& \geq \#\left\{n \in A_{N_k}: \lVert n\alpha \rVert \leq \frac{(\log N_k)^{1+\varepsilon}}{N_k}\frac{N_k}{N}\log \log a_{N_k} \right\} \\
& \geq \#\left\{n \in  A_{N_k}: \lVert n\alpha \rVert' \leq \frac{(\log N_k)^{1+\varepsilon}}{N_k}\frac{N_k}{N}\log \log a_{N_k} \right\}
\\&\geq \mathcal{N}'(N_k,\alpha) \geq \frac{\Psi_{N_k}'}{2} \geq 1,
\end{split}
\]
which proves the statement.

\subsection{Proof of Corollary \ref{cor_alln}}
 
The proof follows exactly the same lines as the proof of Corollary \ref{unif_thm_cases}. Indeed, it is sufficient to show that, for almost every $\alpha \in \mathbb{R}$, for all sufficiently large $N \in \mathbb{N}$, we have
\begin{equation}\label{upp_delta_log}\delta_{\min}^{\alpha}(N) \leq \frac{(\log N)^{2+\varepsilon}}{C_N},\end{equation}
as well as 
\begin{equation}\label{upp_delta_loglog} \delta_{\min}^{\alpha}(N) \leq \frac{(\log N)^{1+\varepsilon}}{C_N} \log \log a_N.\end{equation}
For a set $\A$ we denote by $\mathcal{A}^+$ the set of all positive elements of $\A$ and $C_N^+:= \frac{C_N-1}{2}$ the cardinality of $(\A-\A)^+$.
We define the sequence $(c_n)_{n \in \N}$ by 
$(c_n)_{1 \leq n \leq C_N^+} = (\A_N-\A_N)^+$, while we arrange
$(c_n)_{C_{N-1}^+ < n \leq C_N^+}$ in an arbitrary but fixed way (say, increasing).
We now apply Theorem \ref{main_thm_unif} to that new sequence, and obtain for all $M$ sufficiently large, 
\[\lVert c_n \alpha \rVert \leq \frac{(\log M)^{2+\varepsilon}}{10 M}, \quad 1 \leq n \leq M.\]
Thus in particular, for $M = C_N^+$, and observing that $C_N^+ \leq N^2$, this implies for $N$ sufficiently large,
\[\lVert c_n \alpha \rVert \leq \frac{(\log C_N^+)^{2+\varepsilon}}{10 C_N^+}
\leq \frac{3^{2+\varepsilon}}{10} \frac{(\log N)^{2+\varepsilon}}{C_N}
, \quad 1 \leq n \leq C_N^+,\]
proving \eqref{upp_delta_log}. For \eqref{upp_delta_loglog}, we observe that since $(a_n)_n$ is monotonically increasing,
$\max_{n \leq C_N^+}c_n \leq a_N$, hence for all $n \leq C_N^+$, $\frac{c_n}{\varphi(c_n)} \ll \log \log a_N$. The proof now follows verbatim the proof of Corollary \ref{unif_thm_cases}, with $(c_n)_n$ in the role of $(a_n)_n$.

\bibliographystyle{abbrv}
\bibliography{bibliography}

\Addresses

\end{document}